\documentclass[12pt,reqno]{amsart}
\usepackage{amsmath,amsthm,amsfonts,amssymb,amscd}
\usepackage[latin1]{inputenc}
\usepackage{psfrag,pslatex}
\usepackage{epsfig}
\usepackage{a4wide}
\usepackage{color}

\numberwithin{equation}{section}

\newcommand{\dist}{\operatorname{dist}}

\newcommand{\inter}{\operatorname{int}}
\newcommand{\clos}{\operatorname{closure}}
\newcommand{\diff}{\operatorname{Diff}}
\newcommand{\emb}{\operatorname{Emb}}
\newcommand{\graph}{\operatorname{Graph}}
\newcommand{\diam}{\operatorname{diam}}
\newcommand{\supp}{\operatorname{supp}}

\newcommand{\de} {\delta}       \newcommand{\De}{\Delta}
\newcommand{\ep} {\varepsilon}

\newcommand{\ze} {\zeta}
\renewcommand{\th} {\theta}

\newcommand{\vfi}{\varphi}

\newcommand{\ups}{\upsilon}

\def \NN {{\mathbb N}}
\def \SS {{\mathbb S}}

\def \EE {{\mathbb E}}
\def \FF {{\mathbb F}}
\def \ZZ {{\mathbb Z}}
\def \TT {{\mathbb T}}
\def \RR {{\mathbb R}}
\def \BB {{\mathbb B}}

\newcommand{\qand}{\quad\text{and}\quad}

\def \cB {{\mathcal B}}

\def \cA {{\mathcal A}}
\def \cC {{\mathcal C}}
\def \cH {{\mathcal H}}
\def \cG {{\mathcal G}}

\newcommand{\cV}{{\mathcal V}}
\newcommand{\cW}{{\mathcal W}}

\newcommand{\cU}{{\mathcal U}}
\newcommand{\cP}{{\mathcal P}}

\newcommand{\cM}{{\mathcal M}}
\newcommand{\cD}{{\mathcal D}}
\newcommand{\cK}{{\mathcal K}}

\newtheorem{maintheorem}{Theorem}

\newtheorem{theorem}{Theorem}[section]

\newtheorem{proposition}[theorem]{Proposition}

\newtheorem{lemma}[theorem]{Lemma}

\newtheorem{example}{Example}

\theoremstyle{remark}
\newtheorem{remark}{Remark}

\begin{document}

\author{Vitor Araújo}

\address{Centro de Matemática da
  Universidade do Porto, Rua do Campo Alegre 687, 4169-007
  Porto, Portugal --- Presently at: Instituto de Matematica,
Universidade Federal do Rio de Janeiro,
C.P. 68.530, CEP 21.945-970,
Rio de Janeiro, R. J. , Brazil} 
\email{vdaraujo@fc.up.pt \text{or} vitor.araujo@im.ufrj.br}
\urladdr{http://www.fc.up.pt/cmup/vdaraujo}

\thanks{Work (partially) supported by the Centro de
  Matemática da Universidade do Porto (CMUP), financed by
  FCT (Portugal) through the programmes POCTI (Programa
  Operacional "Ciência, Tecnologia, Inovação") and POSI
  (Programa Operacional Sociedade da Informação), with
  national and European Community structural funds.  V.A.
  was also partially supported by grant FCT/SAPIENS/36581/99
  and enjoyed a post-doc period at PUC (Rio de Janeiro)
  during the preparation of this
  work.\\
  A.T was partially supported by FAPESP FAPESP-Proj.
  Tematico 03/03107-9 and CNPq (Projeto Universal). }

\author{Ali Tahzibi}

\address{Departamento de Matem\'atica,
  ICMC-USP São Carlos, Caixa Postal 668, 13560-970 São
  Carlos-SP, Brazil.}
\email{tahzibi@icmc.sc.usp.br}
\urladdr{http://www.icmc.sc.usp.br/$\sim$tahzibi}



\renewcommand{\subjclassname}{\textup{2000} Mathematics Subject Classification}

\date{\today}

\setcounter{tocdepth}{2}

\title{Physical measures at the boundary of hyperbolic maps}

\begin{abstract}
  We consider diffeomorphisms of a compact manifold with a
  dominated splitting which is hyperbolic except for a
  "small" subset of points (Hausdorff dimension smaller than one, e.g. a denumerable subset) and prove the existence of
  physical measures and their stochastic stability.  The
  physical measures are obtained as zero-noise limits which
  are shown to satisfy the Entropy Formula.
\end{abstract}

\maketitle


\section{Introduction}

Let $M$ be a compact and connected Riemannian manifold  and
$\diff^{1+\alpha}(M)$ be the space of $C^{1+\alpha}$
diffeomorphisms of $M$ for a fixed $\alpha>0$. We write $m$
for some fixed measure induced by a normalized volume form
on $M$ that we call \emph{Lebesgue measure}, $\dist$ for the
Riemannian distance on $M$ and $\|\cdot\|$ for the induced
Riemannian norm on $TM$.

 We say that an invariant probability measure $\mu$ for a
 transformation $f_0:M\to M$ on a manifold $M$ is
 \emph{physical} if the \emph{ergodic basin}
\[
B(\mu)=\left\{x\in M:
\frac1n\sum_{j=0}^{n-1}\varphi(f_0^j(x))\to\int\varphi\,
d\mu\mbox{  for all continuous  } \varphi: M\to\RR\right\}
\]
has positive Lebesgue measure. These measures describe the
asymptotic average behavior of a large subset of points of
the ambient space and are the basis of the understanding of
dynamics in a statistical sense.
It is a challenging problem in the Ergodic Theory of
Dynamical Systems to prove the existence of such invariant
measures.

Let $T_{\Omega(f_0)} = E^s \oplus E^u$ be a hyperbolic
$Df_0$-invariant decomposition (Whitney sum) of the tangent
bundle of the non-wandering set $\Omega(f_0)$ of $f_0$.  The
classical construction of physical measures involves
$f_0$-invariant measures which are absolutely continuous
with respect to Lebesgue measure along the unstable
direction through the points of $\Omega(f_0)$. These
uniformly hyperbolic dynamical systems were the first
general class of systems where these measures were shown to
exist \cite{BR75,Ru76,Si72}.

An invariant probability measure is called \emph{SRB
(Sinai-Ruelle-Bowen) measure}, if it admits positive Lyapunov
exponents and its conditional measures along the unstable
manifolds (in the sense of Pesin theory~\cite{Pe76,FHY83}) are absolutely
continuous with respect to  Lebesgue measure induced on the
unstable manifolds. For a class of dynamical systems
which includes uniformly hyperbolic systems the notions of
\emph{physical} and \emph{SRB} measures coincide.

The SRB measures as defined above are related to a class of
equilibrium states of a certain potential function. Let
$\phi : M \rightarrow \mathbb{R}$ be a continuous function.
Then a $f_0$-invariant probability measure $\mu$ is a
\emph{equilibrium state for the potential $\phi$} if
\[
 h_{\mu}(f_0) + \int \phi \, d\mu =
\sup_{\nu \in \cM} \left\{ h_{\nu}(f_0) + \int \phi \, d\nu \right\},
\]
where $\cM$ is the set of all $f_0$-invariant probability measures.

For uniformly hyperbolic diffeomorphisms it turns out that
physical (or SRB) measures are the equilibrium states for
the potential function $\phi (x) = - \log | \det Df | E^u
(x)|$.  It is a remarkable fact that \emph{for uniformly
  hyperbolic systems these three classes of measures
  (physical, SRB and equilibrium states) coincide}.

We will address the problem of the existence of physical
measures on the boundary of uniformly hyperbolic
diffeomorphisms. The idea is to add small random noise to a
deterministic system $f_0$ in the boundary of uniformly
hyperbolic systems and, for a large class of such maps, we
prove that as the level of noise converges to zero, the
stationary measures of the random system tend to equilibrium
states for $f_0$ which are physical measures. The stationary
measures exist in a very general setting, but the ``zero
noise'' limit measures are not necessarily physical
measures. The specific choice of random perturbation is
important to obtain 
equilibrium states as zero noise limits. These equilibrium
states satisfy Pesin's Entropy Formula and by the
characterization of measures satisfying this formula (whose
proof in~\cite{LY85} is well known to be valid also for
$C^{1+\alpha}$ diffeomorphisms) we deduce that such zero
noise limits are \emph{SRB} measures. In the setting of the
main theorems every $SRB$ measure is a physical measure.
The same general idea has been used in~\cite{CoYo2004} to
obtain \emph{SRB} measures for partially hyperbolic maps
under strong asymptotic growth conditions on every point.

Let $(\theta_\ep)_{\ep>0}$ be a family of Borel probability
measures on $(\diff^{1+\alpha}(M),
\cB(\diff^{1+\alpha}(M)))$, where we write $\cB(X)$ the
Borel $\sigma-$algebra of a topological space $X$. We will
consider random dynamical systems generated by independent
and identically distributed diffeomorphisms with
$\theta_\ep$ the probability distribution driving the choice
of the maps.

We say that a probability measure $\mu^\ep$ on $M$ is
\emph{stationary for the random system $(\diff^{1+\alpha}(M),
  \theta_\ep)$} if
\begin{equation}
  \label{eq:1}
  \int\!\int \varphi(f(x)) \,d\mu^\ep(x)d\theta_{\ep}(f) =
  \int \varphi \,d\mu^\ep \quad \mbox{for all
  continuous  } \varphi:M\to\RR.
\end{equation}

We will assume that $\supp(\theta_\ep)\to f_0$ when $\ep \to
0$ in a suitable topology.  A result based on classical
Markov Chain Theory (see~\cite{Ki88} or~\cite{Ar00}) ensures
that \emph{every weak$^*$ accumulation point of stationary
  measures $(\mu^\ep)_{\ep>0}$ when $\ep \to 0$ is a
  $f_0$-invariant probability measure,} called a
\emph{zero-noise limit measure.}  It is then natural to
study the kind of zero noise limits that can arise and to
define stochastic stability when the limit map $f_0$ admits
physical measures.

We say that a map $f_0$ is \emph{stochastically stable}
(under the random perturbation given by
$(\theta_\ep)_{\ep>0}$) if \emph{every accumulation point $\mu$ of
the family of stationary measures $(\mu^\ep)_{\ep>0}$, when
$\ep\to0$, is a linear convex combination of the physical
measures of $f_0$}.

Stochastic stability has been proved for uniformly expanding
maps and uniformly hyperbolic
systems~\cite{Ki86a,Ki88,Vi97b,Yo85}.  For some
non-uniformly hyperbolic systems, like quadratic maps, H\'enon
maps and Viana maps, stochastic stability has been obtained
much more recently~\cite{AA03,BaV96,BeV2}. The authors have
studied random perturbations of \emph{intermittent maps} and
have proved stochastic stability for these maps for some
parameters and certain types of random
perturbations~\cite{ArTah}. The techniques used were
extended to higher dimensional local diffeomorphisms
exhibiting expansion except at finitely many points,
enabling us to obtain physical measures directly as
zero-noise limits of stationary measures for certain types
of random perturbations, proving also the stochastic
stability of these measures.

Stochastic stability results for maps of the 2-torus which
are essentially Anosov except at finitely many points were
obtained in~\cite{CoYo2004}, the physical probability
measures of which were constructed in a series of papers
using different
techniques~\cite{hu-young1995,hu2001,hu2000}.
Similar results for different kinds of bifurcations away
from Anosov maps at fixed points were also studied in
\cite{CatEnr2001}.

Using ideas akin to~\cite{ArTah} and~\cite{CoYo2004} we
prove the existence of physical probability measures and
their stochastic stability for diffeomorphisms which are
``almost Anosov'' under some geometric and dynamical
conditions.

\subsection{Statement of the results}
\label{sec:statement-results}

We assume that $f_0:U_0\to f_0(U_0)$ is a $C^{1+\alpha}$
diffeomorphism on a relatively compact open subset $U_0$ of a
manifold $M$ which is strictly invariant, that is,
$\clos(f(U_0))\subset U_0$.
During the rest of this paper we set
  $\Lambda=\cap_{n\ge0}\clos{f_0^n(U_0)}$.
Moreover we suppose there exists a
continuous dominated splitting $E\oplus F$ of $T_{U_0} M$
which is $Df_0$-invariant over $\Lambda$,
i.e., there exists $\lambda_0\in(0,1)$ such that for all $x\in
U_0$
\begin{equation}
  \label{eq:domination}
\|Df\mid E(x)\|\cdot \|(Df\mid F(x))^{-1}\|\le\lambda_0.
\end{equation}
We may see $E\oplus F$ on $U_0$ as a continuous extension of
$E\oplus F$ on $\Lambda$.
This assumption ensures the existence of
\begin{description}
\item[stable cones] $\EE^a_x=\{(u,v)\in E(x)\oplus
  F(x): \|v\|\le a\cdot\|u\|\}$;
\item[unstable cones] $\FF^b_x=\{(u,v)\in E(x)\oplus
  F(x): \|u\|\le b\cdot\|v\|\}$;
\end{description}
for all $x\in U$ and $a,b\in(0,1)$, which are $Df_0$-invariant in
the following sense
\begin{itemize}
\item if $x,f_0^{-1}(x)\in U$, then $Df_0^{-1}(\EE^{ a}_x)\subset
  \EE^{\lambda_0 a}_{f_0^{-1}(x)}$;
\item if $x,f_0(x)\in U$, then $Df_0(\FF^b_x)\subset
  \FF^{\lambda_0 b}_{f_0(x)}$;
\end{itemize}
Continuity enables us to unambiguously denote $d_E=\dim (E)$
and $d_F=\dim (F)$, so that $d=d_E+d_F=\dim(M)$. Domination
guarantees the absence of tangencies between stable and
unstable manifolds, since the angles between the $E$ and $F$
directions are bounded from below away from zero at every
point.
Let us fix the unit balls of dimensions $d_E,d_F$
\[
\BB_E=\{w\in\RR^{d_E}: \|w\|_2\le1\}\qand
\BB_F=\{w\in\RR^{d_F}: \|w\|_2\le1\}
\]
where $\|\cdot\|_2$ is the standard Euclidean norm on the
corresponding Euclidean space.
We say that a $C^{1+\alpha}$ embedding $\De:\BB_E\to M$
(respectively $\De:\BB_F\to M$) is a $E$-disk (resp. $F$-disk)
if the image of $D\De(w)$ is  contained in $\EE^a_{\De(w)}$
for all $w\in\BB_E$ (resp. $D\De(w)(\RR^{d_F})\subset
\FF_{\De(w)}^b$ for every $w\in\BB_F$).

In what follows we denote by $\cH(A)$ the
\emph{Hausdorff dimension} of a subset $A\subset M$.
We first state the results without mentioning random
perturbations.

\begin{maintheorem}
  \label{th:limitSRB}
  Let $f_0:M\to M$ be a $C^{1+\alpha}$ diffeomorphism
  admitting a strictly forward invariant open set $U_0$ endowed
  with a dominated splitting $E\oplus F$ such that
\begin{enumerate}
\item $\|Df_0\mid E(x)\|\le1$ and $\|(Df_0\mid
  F(x))^{-1}\|\le1$ for all $x\in U_0$;
\item $F_1=\{x\in U_0: \|(Df_0\mid F(x))^{-1}\|=1\}$
  and $E_1=\{x\in U_0: \|Df_0\mid E(x)\|=1\}$ satisfy
 \[
   \cH(\De\cap E_1)<1  \qand \cH(\hat\De\cap F_1)<1,
 \]
where $\De$ is any $E$-disk and $\hat\De$ is any $F$-disk
contained in $U_0$;

\item $|\det (Df_0\mid F(x))|>1$ for every $x\in F_1$.
\end{enumerate}

If in addition $f_0\mid\Lambda$ is transitive, then there
exists a unique physical measure supported in $\Lambda$,
with $d_F$ positive Lyapunov exponents along the
$F$-direction and whose basin has full Lebesgue measure in
$U_0$.
\end{maintheorem}

We note that if $E_1\cup F_1$ contains no periodic points
and is \emph{finite}, then some power of $f_0$ is a
uniformly hyperbolic map, in which case the conclusions of
Theorem~\ref{th:limitSRB} are known.
Moreover from the dominated decomposition
assumption~\eqref{eq:domination} we easily see that
$E_1\cap F_1=\emptyset.$
We remark also that \emph{the conditions on $E_1$ and $F_1$ in the
statement of Theorem~\ref{th:limitSRB} are automatically
satisfied whenever $E_1$ or $F_1$ is denumerable}.

The restriction on the Hausdorff Dimension is used to show
that any curve inside a $E$-disk (or $F$-disk) intersects
$E_1$ (or $F_1$, respectively) in a zero Lebesgue measure
subset. In particular our results can be obtained assuming
that $E_1$ and $F_1$ do not contain any such curves.

We clearly may specialize this result for a transitive
$C^{1+\alpha}$-diffeomorphism admitting a dominated splitting on the
entire manifold and satisfying items (1)-(3) of
Theorem~\ref{th:limitSRB}, up to replacing $U_0$ and
$\Lambda$ by $M$.

\begin{remark}\label{rmk:partialhyp}
  We can adapt the statement of Theorem~\ref{th:limitSRB} to
  the setting where $U_0$ has a \emph{partially hyperbolic
    splitting}, that is, the strictly forward
  $f_0$-invariant open subset $U_0$ admits a continuous
  splitting $T_{U_0} M=E^s\oplus E^c\oplus E^u$ such that
  \begin{itemize}
  \item both $(E^s\oplus E^c)\oplus E^u$ and $E^s\oplus
    (E^c\oplus E^u)$ are dominated decompositions;
  \item $E^s$ is uniformly contracting and $E^u$ is
    uniformly expanding: there exists $\sigma>1$ satisfying
    $\|Df\mid E^s(x)\|\le\sigma^{-1}$ and $\|(Df\mid
    E^u(x))^{-1}\|\le\sigma^{-1}$ for all $x\in U$;
  \item the restriction of the splitting to $\Lambda$ is
    $Df_0$-invariant.
  \end{itemize}

  If we assume that $\Lambda$ is transitive and either
    \begin{enumerate}
    \item[(1)] $\|Df_0\mid E^c(x)\|\le1$ for all $x\in U_0$,
      and $K=\{x\in U_0: \|Df_0\mid E^c(x)\|=1\}$ is such
      that $\cH(\De\cap K)<1$ for every $E^s\oplus E^c$-disk
      $\De$ contained in $U_0$;
    \end{enumerate}
    or
    \begin{enumerate}
    \item[(2)] $\|(Df_0\mid E^c(x))^{-1}\|\le1$ for all
      $x\in U_0$; $K=\{x\in U_0: \|(Df_0\mid
      E^c(x))^{-1}\|=1\}$ satisfies $\cH(\De\cap K)<1$ for
      every $E^c\oplus E^u$-disk $\De$ contained in $U_0$;
      and $|\det (Df_0\mid E^c(x))|>1$ for every $x\in K$;
    \end{enumerate}
    then there exists a unique absolutely continuous
    $f_0$-invariant probability measure $\mu_0$ with $\dim
    E^u$ (case 1) or $\dim E^c+\dim E^u$ (case 2) positive
    Lyapunov exponents, whose support is contained in
    $\Lambda$ and with an ergodic basin of full Lebesgue
    measure in $U_0$.

    The statement essentially means that \emph{if an
      attractor admits a partially hyperbolic splitting
      which is volume hyperbolic, does not admit mixed
      behavior along the central direction and the neutral
      points along the central direction form a small
      subset, then there exists a physical measure.}
\end{remark}

Cowieson-Young \cite{CoYo2004} have obtained similar
results, albeit for the existence of $SRB$ measures and not
necessarily for \emph{physical} ones (see V\'asquez
\cite{vasquez2006} where it is shown that in the setting of
Remark~\ref{rmk:partialhyp} physical measures will
necessarily be $SRB$ measures).  Moreover Cowieson-Young
obtained a strong result of existence of $SRB$ measure for
partially hyperbolic maps with one-dimensional central
direction using the same strategy of proof. Here we get rid
of the dimensional restriction assuming a dynamical
restriction: see Example \ref{ex:ex-center-dim2} in
Section~\ref{sec:exampl-maps-sett} for a partially
hyperbolic map with two dimensional center-stable and
center-unstable directions $E$ and $F$ in the setting of
Theorem~\ref{th:limitSRB}.

These results will be derived from the following more
technical one, but also interesting in itself.

\begin{maintheorem}
\label{th.limitequilibrium}
Let $f_0:M\to M$ be a $C^{1+\alpha}$ diffeomorphism
admitting a strictly forward invariant open set $U_0$ with a
dominated splitting satisfying items (1)-(2) of
Theorem~\ref{th:limitSRB}.
Then
\begin{enumerate}
\item for any non-degenerate isometric random perturbation
$(\th_\ep)_{\ep>0}$ of $f_0$, every weak$^*$ accumulation
point $\mu$ of a sequence $(\mu^\ep)_{\ep>0}$ of stationary
measures of level $\ep$, when
$\ep\to0$, is an equilibrium state for the potential $-\log
|\det( Df_0\mid F(x))|$, i.e.
\begin{equation}
  \label{eq:3}
  h_\mu(f_0) = \int \log |\det (Df_0\mid F(x))| \, d\mu (x).
\end{equation}
\item every equilibrium state $\mu$ as above is a convex
linear combination of
\begin{enumerate}
\item at most finitely many ergodic equilibrium states
  having positive entropy with $d_F$ positive Lyapunov
  exponents, with
\item probability measures having zero entropy whose support
has constant unstable Jacobian equal to one, i.e., measures
whose Lyapunov exponents are non-positive.
\end{enumerate}
\item every equilibrium state with positive entropy is
 a physical measure for $f_0$.
\item if \emph{the attractor
 $\Lambda$ is transitive}, then there exists at most one
 equilibrium state with positive entropy.
\end{enumerate}
\end{maintheorem}

\begin{remark}
\label{rmk:F1countable}
We note that if $F_1$ is denumerable, then necessarily the
measures in item (2b) of Theorem~\ref{th.limitequilibrium}
are Dirac measures concentrated on periodic orbits whose
tangent map has only non-positive eigenvalues.
\end{remark}

The restriction on the random perturbations means the
following. We assume that $U_0\subset M$ admits an open subset
$V\subset\clos(V)\subset U$ and an action
$\cV\to M$, where $\cV$ is a small neighborhood of the
identity $e$ of a locally compact Lie group $G$ such that
for all $x\in \clos(V)$, setting $g_x:\cV\to M, \, v\mapsto
v\cdot x$, we have
\begin{description}
\item[P1] $g_x(\cV)\subset U_0$;
\item[P2] $g_x(W)$ is a neighborhood of $x$ for every open
  subset $W\subset\cV$;
\item[P3] for every fixed $v\in\cV$ the map
  $g_v:V\to U_0,\, x\mapsto v\cdot x$ is an isometry.
\end{description}
Then we define
\begin{equation}
  \label{eq:additive}
\hat f:\cV\times M \to M, \quad (v,x)\mapsto v\cdot f_0(x)
\end{equation}
and take a probability measure $\th_\ep$  on $\cV$, which
translates into a probability measure on the family
$(\hat f_v)_{v\in\cV}$.

This kind of families are a special case of non-degenerate
random perturbations.  For more on non-degenerate random
perturbation and for examples of non-degenerate isometric
random perturbations, see Section~\ref{sec:non-deg-pert}.
In particular, \emph{Theorems~\ref{th:limitSRB}
  and~\ref{th.limitequilibrium} apply to a bounded
  topological attracting set for a diffeomorphism on a
  domain of any Euclidean space}. In
Section~\ref{sec:non-deg-pert} we show that \emph{this is
  enough to obtain Theorems~\ref{th:limitSRB}
  and~\ref{th.limitequilibrium} in full generality} through
a tubular neighborhood construction.  In particular,
Theorem~\ref{th.limitequilibrium} shows that \emph{in the
  setting of Theorem~\ref{th:limitSRB} (or
  Remark~\ref{rmk:partialhyp}) the physical measures
  obtained are stochastically stable}, as explained in
Section~\ref{sec:stochstability}.


This paper is organized as follows. In
Section~\ref{sec:exampl-maps-sett} we present some examples
in the setting of the main theorems. We outline some general
results concerning random maps in
Section~\ref{sec:non-deg-pert}. In
Section~\ref{sec:stable-unst-manif} we derive the main
dynamical consequence of our assumptions and then, in
Section~\ref{sec:EquiStates}, we prove that equilibrium
states for $f_0$ must be either physical measures or
measures with no expansion.  Finally we construct
equilibrium states using zero-noise limits in
Section~\ref{sec:zero-noise-limits} and put together the
results concluding stochastic stability for $f_0$ and
proving Theorem~\ref{th:limitSRB} in
Section~\ref{sec:stochstability}.


\section{Examples of maps in the setting of the main
  theorems}
\label{sec:exampl-maps-sett}

\begin{example}
  \label{ex:Hu}
  Let $f_0:\TT^2\to\TT^2$ be a $C^{2}$ diffeomorphism with
  $\TT=\SS^1\times\SS^1$, obtained from an Anosov linear
  automorphism of the $2$-torus by weakening the expanding
  direction $F$ of the fixed point $p$ in such a way that
  $Df_0(p)\mid F= Id \mid F$. The stable direction $E$
  continues to be uniformly contracting throughout and $F$
  is still expanded by $Df_0$ on $\TT^2\setminus\{p\}$.

  This kind of maps where studied by Hu and Young
  \cite{hu2000,hu2001,hu-young1995}.  
  In this setting the only physical probability measure for
  $f_0$ is $\de_p$, whose basin contains Lebesgue almost
  every point of $\TT$. Hence
  Theorem~\ref{th.limitequilibrium} shows in particular that
  $\de_p$ is stochastically stable.

  The construction can be adapted to provide maps with
  finitely many periodic orbits with neutral behavior along
  the $F$ direction. We note that $E_1=\emptyset$.
\end{example}

\begin{example}
  \label{ex:productHu}
  Let us take the product $f_0\times E_d$, where $f_0$ is
  given by Example~\ref{ex:Hu} and $E_d:\SS^1\to\SS^1,
  x\mapsto d\cdot x \mod \ZZ$, identifying $\SS^1$ with
  $\RR/\ZZ$ and letting $d\in\NN, d\ge2$. Then
  $E_1=\emptyset$, $E^s=E\times\{0\}$, $E^c=F\times\{0\}$
  and, for big enough $d\ge2$, $E^u=\{(0,0)\}\times\RR$.
  Moreover $F_1=\{p\}\times \SS^1$ and
  $W^u_{loc}(p)\times\SS^1$ is a $E^c\oplus E^u$-disk that
  contains $F_1$, and also $\cH(F_1)=1$.

  In this example $\mu=\de_p\times\lambda$ is the unique
  physical measure, has positive entropy and only one positive
  Lyapunov exponent, where $\lambda$ is Lebesgue measure on
  $\SS^1$.
\end{example}

We note that Examples~\ref{ex:Hu} and~\ref{ex:productHu} can
be seen as ``derived from Anosov'' (DA) maps
\cite{williams1970,Ca93} at the boundary of the set of Anosov
diffeomorphisms.

\begin{example}
  \label{ex:intermitentsolenoid}
Let $f_0:\SS^1\times\RR^2\to\SS^1\times\RR^2, (x,\rho
e^{i\theta})\mapsto (g_\alpha(x), (\rho/10+1/2)\cdot
e^{i(\theta+g_\alpha(x))})$ where again $\SS^1=\RR/\ZZ$
and in $\RR^2$ we use polar coordinates. If $g_\alpha:\SS^1\to\SS^1$ is an
expanding map, then we have the standard solenoid map. Here
we take the $C^{1+\alpha}$ map
\[
g_\alpha(x)= \left\{
\begin{array}{ll}
 x + 2^{\alpha} x^{1+\alpha}, \qquad &  x \in [0,\frac{1}{2})\\
 x - 2^{\alpha} (1-x)^{1 + \alpha},  &  x \in [\frac{1}{2},1]
\end{array} \right.
\]
for $0<\alpha<1$. It is known~\cite{thaler1980} that
$g_\alpha$ admits a unique absolutely continuous invariant
probability measure $\mu$. If $\pi:\SS^1\times \RR \to
\SS^1$ is the natural projection, then on the attractor
$\Lambda=\cap_{n\ge1} f_0^n(\SS^1\times \clos(B(0,1)))$
there is a unique measure $\nu$ such that $\pi_*(\nu)=\mu$,
which is physical and whose basin contains Lebesgue almost
every point of $U_0=\SS^1\times \clos(B(0,1))$.

In this case $F_1=\{0\}\times\RR^2$ but every $F$-disk $\De$
intersects $F_1$ at most finitely many times, since $\De$
must be locally a graph over
$S^1$. Theorem~\ref{th:limitSRB} holds and
Theorem~\ref{th.limitequilibrium} shows that every
equilibrium state is a convex linear combination of
$\de_{(0,5/9)}$ with $\nu$ ($(0,5/9)$ is the unique fixed
point of $f_0$).

If we let $\alpha\ge1$, then $g_\alpha$ is of class $C^2$
and~\cite{thaler1983} $\de_0$ is the unique physical measure
for $g_\alpha$. Theorem~\ref{th.limitequilibrium} shows that
$\de_{(0,5/9)}$ is stochastically stable. Since
$\pi_*(\de_{(0,5/9)})=\de_0$ it is not difficult to see that
$\de_{(0,5/9)}$ has basin containing $U_0$ Lebesgue modulo
zero, so $\de_{(0,5/9)}$ is the physical measure for
$\Lambda$.
\end{example}

\begin{example}
\label{ex:produtotorto}
Let $f_0:\TT\times\RR^2\to\TT\times\RR^2, (t,x,\rho
e^{i\theta})\mapsto (E_d(t), g(E_d(t),x),(\rho/10+1/2)\cdot
\exp[i(\theta+g(E_d(t),x))])$, where $E_d$ was defined in
Example~\ref{ex:productHu}, $d\in\NN, d\ge2$ and
$g:\TT\to\SS^1$ of class $C^{1+\alpha}$ is an extension of
$g_\alpha$ from Example~\ref{ex:intermitentsolenoid} to
$\TT$ given by
\[
  g(t,x)=\left\{
\begin{array}{ll}
x(1+0.1\cdot\sin^2(\pi t)) + 2^{\alpha}(1-0.1\cdot\sin^2(\pi t)) x^{1+\alpha},
\qquad
&  x \in [0,\frac{1}{2})
\\
1-(1-x)(1+0.1\cdot\sin^2(\pi t))-
2^{\alpha}(1-0.1\cdot\sin^2(\pi t))(1-x)^{1 + \alpha},
&  x \in [\frac{1}{2},1]
\end{array} \right.
\]
for some fixed $0<\alpha<1$, where $\SS^1=\RR/\ZZ$. Then we
have $E^s=\{(0,0)\}\times\RR^2$,
$E^c=\{0\}\times\SS^1\times\{0\}$ and
$E^u=\SS^1\times\{0\}\times\{0\}$ for big enough $d$. The
conditions on item 2 of Remark~\ref{rmk:partialhyp} hold
with $K=E_d^{-1}(\{0\})\times\{0\}\times\RR^2$, because
$|g_t'(x)|=|D_2 g(t,x)|\ge1$ and equals $1$ only at $(0,0)$.

The natural projection $\pi:\TT\times\RR^2\to\TT$ conjugates
$f_0$ to $f_1:\TT\to\TT, (t,x)\mapsto (E_d(t),g(E_d(t),x))$
over the attractor $\Lambda=\cap_{n\ge1}f_0^n(\TT\times
\clos(B(0,1)))$. We note that each $g_t$ is conjugate to $E_2$
through a homeomorphism $h_t$ which depends continuously on
$t\in\SS^1$ in the $C^0$ topology. Hence $H(t,x)=(t,
h_t(x))$ is a homeomorphism of $\TT$ such that $H\circ
f_1=(E_d\times E_2)\circ H$ and since $E_d\times E_2$ is
transitive, then $f_1$ and also $f_0$ are transitive.

This shows that we can apply Remark~\ref{rmk:partialhyp}
obtaining the existence of a unique physical measure for
$f_0$.
\end{example}

\begin{example}
  \label{ex:Devil}
Let $\cK \subset I=[0, 1]$ be the middle third Cantor set and
$(a_i^n, b_i^n), i= 1, \cdots , 2^{n-1}$ be an enumeration
for the gaps of the $n$-th generation in the construction of $\cK$.
We define $\beta$ on  any given gap interval $(a, b)$ as
\[
\beta (x)= \left\{
\begin{array}{ll}
 x - a,  & \mbox{if  } x \in (a, \frac{a + b}{2})\\
 b - x,   & \mbox{if  } x \in (\frac{a + b}{2}, b)
\end{array} \right..
\]
Then the map $\beta:I\setminus\cK\to I$ is uniformly
continuous and so  we can continuously extend it to $I$
setting $\beta\mid\cK\equiv 0$. Moreover it is easy to see
that $\beta\mid (I \setminus \cK)$ is Lipschitz (with
Lipschitz constant $1$) and so is its extension to $I$.

It addition, with respect to Lebesgue measure on $I$, we get $\int_0^1
\beta < \infty$ and if $g_0: I \to\RR$ is given by $g_0(x) = x +
\frac{\int_0^x \beta} {\int_0^1 \beta}$, then $g_0 (0) = 0, g_0(1) =
2$ and $g_0$ induces a $C^1$ map of the circle onto itself whose
derivative is Lipschitz satisfying $g_0^\prime\mid\cK\equiv1$ and
$g_0^\prime\mid(I\setminus\cK)>1$.

The map $g_0:\SS^1\to\SS^1$ is mixing since $\sigma(J)=|g_0(J)|/|J|>1$
for every arc $J\subset\SS^1$, where $|\cdot|$ denotes length. Indeed
the continuity of the map $\sigma$ on arcs together with the
compactness of the family $\Gamma(\ell)=\{J\subset\SS^1: J \mbox{ is
an arc and } |J|\ge \ell \}$, for any given bound $\ell>0$ on the
length, show that there exists $\sigma(\ell)>1$ such that
$|g_0(J)|\ge\sigma(\ell)\cdot |J|$ for any given arc $\emptyset\neq
J\subset\SS^1$. Hence for every nonempty arc $J$ there exists
$n=n(J)\in\NN$ such that $g_0^n(J)=\SS^1$.

Replacing $g_\alpha$ by $g_0$ in the definition of $f_0$
within Example~\ref{ex:intermitentsolenoid}, we get a
$C^{1+1}$ map from the solid torus into itself whose
topological attractor satisfies the conditions of Theorem A,
where $F_1$ is Cantor set.

\end{example}

\begin{example}
  \label{ex:ex-center-dim2}
  We present an example of a transitive diffeomorphism with
  2-dimensional center-unstable and center-stable directions
  in the setting of Theorem~\ref{th:limitSRB}. The idea for
  the construction of this example comes from the
  construction of stably transitive diffeomorphisms without
  any uniformly hyperbolic direction in \cite{BoV00}.
\end{example}

We start with a linear Anosov diffeomorphism $f_0$ induced
in $\mathbb{T}^4$ by a linear map of $\mathbb{R}^4$ with
eigenvalues
\[
0 < \lambda_1 < \lambda_2 < \frac{1}{3} < 3 <
\lambda_3 < \lambda_4.
\]
Up to replacing it by some iterate we may suppose that $f_0$
has at least two fixed points $p$ and $q$. For small
$\alpha>0$ we consider a new diffeomorphism $f$ satisfying
the following properties:
 \begin{enumerate}
 \item $f$ has center-unstable cone field $C^{cu}$ and
   center-stable cone field $ C^{cs}$ with width bounded by
   $\alpha > 0$, respectively, containing the unstable and
   stable subbundle of $f_0$;
 \item there exists $\sigma > 1$ such that $|\det Df| TD^{cu}| \geq
 \sigma$ for every disk tangent to the cone field $C^{cu}$ and $|\det Df|
 TD^{cs}| \leq \sigma^{-1}$ for every disk tangent to the cone field
 $C^{cs};$
\item there exist $\lambda \leq 1/3$ such that $\| Df(x)
  v^{cu}\| \geq \lambda^{-1} \|v^{cu}\|$ and $\| Df^{-1}(x)
  v^{cs}\| \geq \lambda^{-1} \|v^{cs}\|$ for every $x$
  outside the union of two small balls $V_p$ around $p$ and
  $V_q$ around $q$, and $v^{cu}\in C^{cu}$ and $ v^{cs}\in
  C^{cs};$
\item the stable index (the dimension of uniformly
  contracting subbundle of the tangent space) at $p$ is equal to
  $1$ and the unstable index is equal to $2.$ For $q$ the
  indexes are given just exchanging "stable" by "unstable" in
  the case of $p$;
\item inside the union of the balls mentioned at item 3
  above we have $\| Df(x) v^{cu}\| \geq \|v^{cu}\|$ and $\|
  Df^{-1}(x) v^{cs}\| \geq \|v^{cs}\|$.
\end{enumerate}

To obtain such $f$ we just modify $f_0$ in a small
neighbourhood along the weaker stable direction of $p$ and
the weaker unstable direction of $q.$ So the strong stable
and strong unstable directions are preserved and $f$ is
partially hyperbolic.  Since $f$ is transitive (see
\cite{BoV00}), by the special tangent bundle decompositions
at $p$ and $q$ we conclude that there cannot exist any
two-dimensional invariant sub-bundle with uniformly
hyperbolic behavior (either uniformly expanding or uniformly
contracting).  In this example, $E_1 = \{p\} $ and $F_1 =
\{q\}$, where $E_1$ and $F_1$ are as in Theorem
\ref{th:limitSRB}, and the tangent bundle admits a dominated
decomposition into \emph{four invariant one-dimensional subbundles}
$E^{ss}\oplus G\oplus H \oplus E^{uu}$, both $E=E^{ss}\oplus
G$ and $F=H\oplus E^{uu}$ are two-dimensional and $f$
satisfies all the hypothesis of Theorem~\ref{th:limitSRB}.

\section{Random perturbations}
\label{sec:non-deg-pert}

Let a parameterized family of maps $\hat f: X\to
\diff^{1+\alpha}(M), t\mapsto f_t$ be given, where $X$ is a
connected compact metric space. We identify a sequence
$f_0,f_1,f_2,\dots$ from $\diff^{1+\alpha}(M)$ with a
sequence $\omega_0,\omega_1,\omega_2,\dots$ of parameters in
$X$ and the probability measure $\theta_\ep$ can be assumed
to be supported on $X$.
We set $\Omega=X^\NN$ to be the space of sequences
$\omega=(\omega_i)_{i\ge0}$ with elements in $X$ (here we
assume that $0\in\NN$). Then we define in $\Omega$ the
standard infinite product topology, which makes $\Omega$ a
compact metrizable space.
The standard product probability measure
$\theta^\ep=\theta_\ep^\NN$ makes $(\Omega,\cB,\theta^\ep)$
a probability space. We write $\cB=\cB(\Omega)$ for the
$\sigma$-algebra generated by cylinder sets: the minimal
$\sigma-$algebra containing all sets of the form $\{ \omega
\in \Omega : \omega_0 \in A_0, \omega_1\in A_2, \omega_2 \in
A_2, \cdots , \omega_l \in A_l \}$ for any sequence of Borel
subsets $A_i \subset X, i= 0, \cdots , l$ and $l\ge1$.
We use the following skew-product map
\[
   F : \Omega\times M  \to \Omega\times M, \quad
  (\omega,x) \mapsto (\sigma(\omega),f_{\omega_0}(x))
\]
where $\sigma$ is the left shift on sequences:
$(\sigma(\omega))_n = \omega_{n+1}$ for all $n\ge0$. It is
not difficult to see that $\mu^\ep$ is a stationary measure
for the random system $(\hat f,\theta_\ep)$ (i.e. satisfying
\eqref{eq:1}) if, and only if, $\theta^\ep \times \mu^\ep$
on $\Omega \times M$ is $F$-invariant.
We say that $\mu^\ep$ is ergodic if $\theta^{\ep} \times
\mu^\ep$ is $F$-ergodic.

If we define $\hat \Omega=X^\ZZ$ to be the set of all
bi-infinite sequences $(\omega_i)_{i\in \ZZ}$ of elements of
$X$, then we can define $G$ to be the invertible natural
extension of $F$ to this space:
\[
   G : \hat\Omega\times M  \to \hat\Omega\times M,\quad
(\omega,x) \mapsto (\sigma(\omega),f_{\omega_0}(x)).
\]
This map is invertible and
$G^{-1}(\omega,x)=(\sigma^{-1}(\omega),f^{-1}_{\omega_{-1}}(x))$.
On $\hat\Omega$ we set the natural product topology and the
product $\sigma$-algebra $\hat\cB=\cB(\hat\Omega)$ generated by
cylinder sets as above but now with indexes in $\ZZ$. The
product probability measure $\hat\th^\ep=\th_\ep^\ZZ$ makes
$(\hat\Omega,\hat\cB,\hat\th^\ep)$ a probability space.  We
set the following notation for the natural projections
\[
\pi_M:\Omega\times M\to M,\quad
\hat\pi_M:\hat\Omega\times M\to M, \quad
\hat\pi_\Omega:\hat\Omega\times M\to\hat\Omega,
\quad\mbox{and}\quad
\hat\pi:\hat\Omega\times M\to\Omega\times M.
\]
For $\omega\in\hat\Omega$ and
for $n \in \ZZ$ we define for all $x\in M$
$$
f^n_{\omega} = (\hat\pi_M \circ G^n) (x)=\left\{
\begin{array}[l]{ll}
(f_{\omega_{n-1}} \circ \dots \circ f_{\omega_0}) (x),
& n>0
\\
x,  & n=0
\\
(f_{\omega_{-n}}^{-1} \circ \dots \circ f_{\omega_{-1}}^{-1}) (x),
& n<0
\end{array}
\right..
$$
Given $x \in M$ and $\omega \in \hat\Omega$ the sequence
$(f^n_{\omega}(x))_{n \geq 1}$ is a {\it random orbit} of
$x$. Analogously we set $f^n_\omega=\pi_M\circ F^n$ for
$n\ge0$ and $\omega\in\Omega$.

From now on we assume that the family $(\theta_{\ep})_ {\ep > 0}$ of
probability measures on $X$ is such that their supports
have non-empty interior
and $\supp (\theta_\ep) \rightarrow \{t_0\}$ when $\ep
\rightarrow 0$, where $t_0\in X$ is such that $f_{t_0}=f_0$.

\subsection{Non-degeneracy conditions}
\label{sec:non-degen-cond}

In what follows we write $ f_x^n : \Omega \rightarrow
M$ for the map $\omega\in\Omega\mapsto f^n_\omega(x)$, for
every $n\ge0$.
We say that $(\hat f, \theta_\ep)_{\ep > 0}$ is a
\emph{non-degenerate random perturbation} of $f_0=f_{t_0}$ if, for
every small enough $\ep$, there is $\delta_1=\delta_1(\ep)> 0$ such
that for  all $x \in U$
   \begin{description}
   \item[ND1] $\{f_t(x): t\in \supp(\theta_\ep) \}$ contains
     a ball of radius $\delta_1$ around $f_{t_0}(x)$;
   \item[ND2] $(f_x)_{*} \theta_{\ep} $ is
     absolutely continuous with respect to $m$.
   \end{description}

\begin{remark}\label{rmk.noatoms}
We note that $\th_\ep$ cannot have atoms by condition  ND2 above.
\end{remark}

The following is a finiteness result for non-degenerate
random perturbations.

\begin{theorem}
  \label{thm:randompert}
  Let $(\hat f, \theta_\ep)_{\ep > 0}$ be a non-degenerate
  random perturbation of $f_0$.  Then for each $\ep>0$ there
  are finitely many absolutely continuous ergodic measures
  $\mu^\ep_1,\dots\mu^\ep_{l(\ep)}$, and for each $x\in U$
  there is a $\th^\ep$ mod $0$ partition
  $\Omega_1(x),\dots,\Omega_{l(\ep)}(x)$ of $\Omega$ such
  that for $1\le i\le {l(\ep)}$
 $$
 \mu_i^\ep=\lim_{n\to+\infty} \frac1n \sum_{j=1}^{n}
 \de_{f_{\omega}^j x} \quad \mbox{ for}\quad \omega \in
 \Omega_i(x).
 $$
 Moreover the interior of the supports of the physical measures are
 nonempty and pairwise disjoint.
\end{theorem}

\begin{proof}
See \cite{Ar00},  \cite{Ze03} or \cite{brin-kifer1987}.
\end{proof}

The continuity of the map $F$ is enough to get the
\emph{forward invariance of $\supp(\mu^\ep)$ for any stationary measure
$\mu^\ep$}, i.e.
if $x\in \supp(\mu^\ep)$ then
$f_t(x)\in\supp(\mu^\ep)$ for all
$t\in\supp(\th_\ep)$,
since $\th^\ep\times\mu^\ep$ is $F$-invariant.  By
non-degeneracy condition ND1 $\supp(\mu^\ep)$ contains a ball of
radius $\de_1=\de_1(\ep)$. Moreover defining the
\emph{ergodic basin} of $\mu^\ep$ by
\[
B(\mu^\ep)=\left\{x\in M:
\frac1n\sum_{j=1}^{n}\varphi(f_\omega^j(x))\to\int\varphi\,
d\mu\mbox{  for all } \varphi\in C(M,\RR) \mbox{ and }
\th^\ep\mbox{-a.e. } \omega\in\Omega \right\},
\]
then $m(B(\mu^\ep))>0$, since $\mu^\ep(B(\mu^\ep))=1$ by the
Ergodic Theorem applied to $(F,\th^\ep\times\mu^\ep)$ and
$\mu^\ep\ll m$.

These non-degeneracy conditions are not too restrictive
since we can always construct a non-degenerate random
perturbation of any differentiable map of a compact manifold
of finite dimension, with $X$ the closed ball of radius 1
around the origin of a Euclidean space, see~\cite{Ar00} and
the following subsection.


\subsection{Isometric random perturbations}
\label{sec:isom-rand-pert}

We present below the two main types of families of maps we
will be dealing with, satisfying conditions P1-P3 stated in
Subsection~\ref{sec:statement-results}.

\begin{example}[Global additive perturbations]
  \label{ex:additiveperturbations}
Let $M$ be a homogeneous space, i.e., a compact
connected finite dimensional Lie Group admitting an
invariant Riemannian metric. Fixing a neighborhood $\cU$
of the identity $e\in M$ we can define a map $f:\cU\times M\to
M, (u,x)\mapsto L_u( f_0(x))$, where $L_u(x)=u\cdot x$ is
the left translation associated to $u\in M$. The invariance
of the metric means that left (an also right) translations
are isometries, hence fixing $u\in \cU$ and taking any
$(x,v)\in TM$ we get
\begin{equation}
  \label{eq:invderivative}
\|Df_u(x)\cdot v\|=\|DL_u(f_0(x))
(Df_0(x)\cdot v)\|=\|Df_0(x)\cdot v\|.
\end{equation}
In the particular case of $M=\TT^d$, the $d$-dimensional
torus, we have $f_u(x)=f_0(x)+u$ and this simplest case
suggests the name \emph{additive random perturbations} for
random perturbations defined using families of maps of this
type. It is easy to see that if the probability measure
$\th_\ep$ is absolutely continuous and supported on a open
subset $X$ of $\cU$, then conditions P1, P2 and P3 are met.
\end{example}

\begin{example}[Local additive perturbations]
  \label{ex:localisometric}
  If $M=\RR^d$ and $U_0$ is a bounded open subset of $M$
  strictly invariant under the diffeomorphism $f_0$, i.e.,
  $\clos(f_0(U_0))\subset U_0$, then
  we can define a non-degenerate isometric random
  perturbation setting
  \begin{itemize}
  \item $V=f_0(U_0)$ (so that $\clos(V)=\clos(f_0(U_0))\subset U_0$);
  \item $G\simeq\RR^d$ the group of translations of $\RR^d$;
  \item $\cV$ a small enough neighborhood of the origin in $G$.
  \end{itemize}
  Then for $v\in\cV$ and $x\in V$ we have $f_v(x)=v\cdot
  x=x+v$, with the standard notation for vector addition,
  and clearly $f_v(x)=x+v$ is an isometry and satisfies both
  conditions P1 and P2.
\end{example}

Now we show that we can construct non-degenerate isometric
random perturbations in the setting of
Examples~\ref{ex:additiveperturbations}
and~\ref{ex:localisometric}.
We define the family of maps
$\hat f$ as in~\eqref{eq:additive}.
The local compactness of $G$ gives a Haar measure $\nu$ on
$G$ and the isometry condition ensures that $\dim (G)=d$ and
that $(\hat f_x)_*(\nu\mid\cV) \ll m$.  Hence for every
probability measure $\th_\ep$ given by a probability density
with respect to $\nu$ we have $(\hat f_x)_*\th_\ep \ll m$,
and this gives condition ND2.

Moreover whenever $\supp(\th_\ep)$ has nonempty interior in
$\cV$ then condition P2, together with the compactness of
$\clos(V)$, ensure that there is $\de=\de(\ep)>0$ such that
condition ND1 is satisfied. Thus we get conditions ND1 and
ND2 choosing $\th_\ep$ as a probability density in $\cV$
whose support has nonempty interior, and setting $X=\cV$ for
the definition of $\Omega,\hat\Omega$.


\subsubsection{Isometric perturbations of maps in arbitrary
manifolds}
\label{sec:Isom-pert-maps}

Now we show that for any given map $f_0$ is the setting of
Theorems~\ref{th:limitSRB} or~\ref{th.limitequilibrium}, we
may define a random isometric perturbation of a particular
extension of $f_0$ as in Example~\ref{ex:localisometric},
which is partially hyperbolic.

We may assume without loss that $M$ is a compact sub-manifold
of $\RR^N$ and that $\|\cdot\|$ and $\dist$ are the ones
induced on $M$ by the Euclidean metric of $\RR^N$, by a
result of Nash~\cite{Nash1954,Nash1956} with $N\ge d(3d+11)/2$. Let $W_0$ be an
open \emph{normal tubular neighborhood} of $M$ in $\RR^N$,
that is, there exists $\Phi:W\to W_0, (x,u)\mapsto x+u$ a
($C^\infty$) diffeomorphism from a neighborhood $W$ of the
zero section of the normal bundle $TM^\perp$ of $M$ to
$W_0$. Let also $\pi:W_0\to M$ be the associated projection:
$\pi(w)$ is the closest point to $w$ in $M$ for $w\in W_0$,
so that the line through the pair of points $w,\pi(w)$ is
normal to $M$ at $\pi(w)$, see e.g. \cite{hirsch1976} or
\cite{guillemin-pollack1974}.  Now we define for
$\rho_0\in(0,1)$
\[
F:W\to W,\quad (x,u)\mapsto (f_0(x),\rho_0\cdot u)
\qand
F_0:W_0\to W_0, \quad w\mapsto (\Phi \circ F \circ \Phi^{-1})(w).
\]
Then clearly $F_0$ is a diffeomorphism onto its image,
$\clos{F_0(W_0)}\subset W_0$ and $M=\cap_{n\ge0}
F_0^n(W_0)$. Moreover if $f_0$ admits a dominated splitting
$E\oplus F$ in a strictly forward $f_0$-invariant set
$U_0\subset M$, then $F_0$ has a dominated splitting
$E^s\oplus E\oplus F$ in the strictly forward
$F_0$-invariant set $\hat U_0=\pi^{-1}(U_0)\subset W_0$,
where $E^s(w)$ is normal to $T_w M$ at $w\in M$ and uniformly
contracted by $DF_0$, as long as $\rho_0$ is close enough
to zero.

We can now define a random isometric perturbation of
$F_0$ and obtain Theorem~\ref{th:limitSRB} as a
corollary of Theorem~\ref{th.limitequilibrium}.  For that \emph{it is
enough to prove Theorem~\ref{th.limitequilibrium} for non-degenerate
random isometric perturbations on an strictly invariant open subset of
the Euclidean space}. Then given $f_0$ we construct $F_0$ as explained
above and note that any $F_0$-invariant measure must be concentrated
on $M\subset\hat U_0$, thus the results obtained for $F_0$ are easily
translated for $f_0$.


\subsubsection{The random invariant set}
\label{sec:RandomInvSet}

In this setting, letting $U_0$ denote the strictly forward
$f_0$-invariant set from the statements in
Section~\ref{sec:statement-results} and $U_k=f_0^k(U_0)$ for
a given $k\ge1$, we have that for some $\ep_0>0$ small enough
\[
\cW=\bigcap_{n\ge0}\clos{G^n(\hat\Omega\times U_k)}
\subset  \hat\Omega\times U_{k-1}
\quad\mbox{and}\quad
\hat\Lambda=\hat\pi_M(\cW)\subset U_{k-1}.
\]
Moreover $\cW$ is $G$-invariant (and $\hat\pi(\cW)$ is
$F$-invariant), where we set $X=\clos{B(0,\ep_0)}$ for
the definition of $\hat\Omega$ (and of $\Omega$).

Indeed we have $\clos{U_k}\subset U_{k-1}$ and
$d_k=\dist(\clos{U_k},M\setminus U_{k-1})>0$. Then we may
find $\ep_0>0$ such that $\dist(f_v(x),f_0(x))\le d_k/4$
for all $v\in B(0,\ep_0)$ and $x\in U_k$.
Hence
\[
f_v(U_k)\subset B\left(\clos(U_k),\frac{d_k}2\right)\subset U_{k-1}
\quad\mbox{for all}\quad v\in B(0,\ep_0),
\]
where $B(A,\de)=\cup_{z\in A}B(x,\de)$ is the
$\de$-neighborhood of a subset $A$, for  $\de>0$.
In addition the $G$-invariance of $\hat\Lambda$ ensures that
\begin{equation}
  \label{eq:invariancianegativa}
  \mbox{if}\quad (\omega,x)\in \cW\quad\mbox{then}\quad
  f_\omega^n(x)\in\hat\Lambda\quad\mbox{for all}\quad n\in\ZZ.
\end{equation}


\subsection{Metric entropy for random perturbations}
\label{sec:metr-entr-rand}

We outline some definitions of metric entropy for random
dynamical systems which we will use and relate them.
Let $\mu^\ep$ be a stationary measure for the random system
given by $(\hat f,\theta_\ep)_{\ep>0}$. Since we are dealing with
randomly chosen invertible maps the following results relating
$F$- and $G$-invariant measures  will be needed.

\begin{lemma}{\cite[Prop. I.1.2]{LQ95}}
  \label{le:chapeumu}
  Every stationary probability measure $\mu^\ep$ of the
  random system given by $(\hat f,\theta_\ep)_{\ep>0}$
  admits a unique probability measure $\hat\mu^\ep$ on
  $\hat\Omega\times M$ which is $G$-invariant and
  $\hat\pi_*(\hat\mu^\ep)=\th^\ep\times\mu^\ep$. Moreover
  $(\hat\pi_\Omega)_*\hat\mu^\ep=\hat\th^\ep$,
  $(\hat\pi_M)_*\hat\mu^\ep=\mu^\ep$ and
  $G^n_*(\hat\th^\ep\times \mu^\ep)$ tends to $\hat\mu^\ep$
  weakly$^*$ when $n\to+\infty$.
\end{lemma}

We will need to consider weak$^*$ accumulation points of
$G$-invariant measures in the following sections, so we
state the following property whose proof follows standard
lines.

\begin{lemma}
  \label{le:weakcontinuous}
  Let $\mu^0$ be a weak$^*$ limit of $\mu^{\ep_k}$ for a
  sequence $\ep_k\to0^+$ when $k\to\infty$. Let $\hat\mu^0$
  be a weak$^*$ accumulation point of the sequence
  $\hat\mu^{\ep_k}$. Then
  $\hat\mu^0=\de_{\omega_0}\times\mu^0$, where
  $\de_{\omega_0}$ is the Dirac mass at
  $\omega_0=(\dots,t_0,t_0,t_0,\dots)\in\hat\Omega$.
\end{lemma}

Here is one possibility of the calculation of the metric
entropy.

\begin{theorem}{\cite[Thm. 1.3]{Ki86}} \label{thm.metr-entr-rand}
 For any finite measurable partition $\xi$ of $M$
 $$
 h_{\mu^\ep}((\hat f,\theta_\ep), \xi) = \inf_{n\ge1}
 \frac{1}{n} \int H_{\mu^\ep} \big(
 \bigvee_{i=-n}^n f^i_{\omega} (\xi) \big) d \theta^{\ep} (\omega)
 $$
 is finite and is called \emph{the entropy of the random
 dynamical system} with respect to $\xi$ and to $\mu^\ep$.
\end{theorem}

We define $h_{\mu^\ep} (\hat f,\theta_\ep)= \sup_\xi \,
h_{\mu^\ep}((\hat f,\theta_\ep), \xi)$ as the
\emph{metric entropy} of the random dynamical system $(\hat
f,\theta_\ep)$, where the supremum is taken over all
measurable partitions.

Let $\cB \times M$ be the minimal $\sigma-$algebra
containing all products of the form $A \times M$ with $A \in
\cB.$ We write $\hat\cB\times M$ for the analogous
$\sigma$-algebra with $\hat\cB$ in the place of $\cB$. We
denote by $h_{ \theta^{\ep}\times \mu^\ep}^{\cB \times M}
(F)$ the conditional metric entropy of the transformation
$F$ with respect to the $\sigma$-algebra $\cB \times M.$
(See e.g.~\cite[Chpt. 0]{LQ95} for a definition and
properties of conditional entropy.) Again we also denote by
$h_{\hat\mu^\ep}^{\hat\cB\times M}(G)$ the conditional
entropy of $G$ with measure $\hat\mu^\ep$ with respect to
$\hat\cB\times M$.

\begin{theorem}{\cite[Prop. I.2.1 \& Thm. I.2.3]{LQ95}}
\label{thm.randentropyFG}
Let $\mu^\ep$ be a stationary probability measure for the
random system given by $(\hat f,\theta_\ep)$. Then
$h_{\mu^\ep}(\hat f,\theta_\ep) = h_{ \theta^{\ep} \times
  \mu^\ep}^{ \cB \times M}
(F)=h_{\hat\mu^\ep}^{\hat\cB\times M}(G)$.
\end{theorem}

The analogous Kolmogorov-Sinai result about generating
partitions is also available in this setting. We let
$\cA=\cB(M)$ be the Borel $\sigma$-algebra of $M$. We say
that a finite partition $\xi$ of $M$ is a \emph{random
  generating partition} for $\cA$ if
$\vee_{i=-\infty}^{+\infty} f_{\omega}^i (\xi) =\cA$
for $\hat\theta^{\ep}$ almost all $\omega\in \hat\Omega$.

\begin{theorem}{\cite[Cor. 1.2]{Ki86}}
  \label{thm.KSrandom}
 Let $\xi$ be a random generating partition for
  $\cA$. Then $h_{\mu^\ep}(\hat
  f,\theta_\ep)=h_{\hat\mu^\ep}^{\hat\cB\times
    M}(G,\hat\Omega\times\xi)$.
\end{theorem}

We note that in \cite{Ki86} this result is stated only for
one-sided sequences. However we know that the
Kolmogorov-Sinai Theorem applied to an invertible
transformation like $G$ demands that a partition $\zeta$ of
$\hat\Omega\times M$ be generating in the sense that
$\vee_{i\in\ZZ} G^i(\zeta)$ equals $\hat\cB\times M,\,
\hat\mu^\ep\bmod 0$. Since we are calculating a conditional
entropy, it is enough that $(\vee_{i\in\ZZ}
G^i(\zeta))\vee(\hat\cB\times M)$ be the trivial partition
in order that $h_{\hat\mu^\ep}^{\hat\cB\times
  M}(G,\zeta)=h_{\hat\mu^\ep}^{\hat\cB\times M}(G)$.  In
particular, for $\zeta=\hat\Omega\times\xi$, we have
$G^i(\zeta)= \{\{\sigma^k(\omega)\}\times
f_\omega^k(\xi): \omega\in\hat\Omega\}$ for $i\in\ZZ$ so
\[
\bigvee_{i=-n}^{n} G^i(\zeta)=\left\{
\{\omega\}\times\bigvee_{i=-n}^n f_\omega^i(\xi) :
\omega\in\hat\Omega \right\}.
\]
Hence $\zeta$ generates $(\hat\Omega\times\cA,\hat\mu^\ep)$
if, and only if, $\xi$ is generating for $\cA$, since
$(\hat\pi_M)_*\hat\mu^\ep=\mu^\ep$.



\section{Expanding and contracting disks}
\label{sec:stable-unst-manif}

Here we derive the main local dynamical properties of the maps in
the setting of Theorem~\ref{th.limitequilibrium}.
We show that $F$-disks  (respectively $E$-disks)  are expanded
(resp. contracted) by the action of $f_0$, and that the
rate of expansion (resp. contraction) is uniform for all
isometrically perturbed $g$ in a $C^{1+\alpha}$-neighborhood
of $f_0$, but depends on the size of the disks.
We also show that the curvature of such disks
remains bounded under iteration. These are consequences of
the domination condition~\eqref{eq:domination} on the
splitting together with non-mixing of expanding/contracting
behavior along the $E$ and $F$ directions given by condition
(1) in Theorem~\ref{th:limitSRB}.

We note that for $g$ sufficiently $C^1$-close to $f_0$ and
for a small $\zeta\in(0,\alpha)$ and a slightly bigger
$\tilde\lambda_0\in(\lambda_0,1)$ we still have for all $x\in U$
\begin{equation}
  \label{eq:dominationclose}
  \|Dg\mid E(x)\|\cdot \|(Dg\mid F(x))^{-1}\|^{1+\zeta}\le\tilde\lambda_0.
\end{equation}
Moreover since $\clos(f_0(U))\subset U$, then for $g$
sufficiently $C^0$-close to $f_0$ in $\diff^{1+\alpha}(M)$
we also have $\clos(g(U))\subset U$, see
Subsection~\ref{sec:RandomInvSet} for more details. We
denote by $\cU$ a $C^{1+\alpha}$ neighborhood of $f_0$ where
all of the above is valid for $g\in\cU$.


\subsection{Domination and bounds on expansion/contraction}
\label{sec:expand-contr-cone}

The domination condition~\eqref{eq:domination} ensures that
the splitting $E(x)\oplus F(x)$ varies continuously in
$\Lambda$ and that there are stable and unstable cone fields
$\EE^a, \FF^b$, already defined in
Subsection~\ref{sec:statement-results} for small $a,b>0$,
which are $Dg$-invariant
for every $g$ sufficiently $C^1$-close to $f_0$. This is a
general property of dominated splittings.

We define a norm on $T_{\clos(U)} M$ more adapted to our purposes
using the splitting: for every $x\in U$ and $w\in T_x M$ we
write
\begin{equation}
  \label{eq:defnorma}
  w=(u,v)\in E(x)\oplus F(x) \mbox{  and set  }
  |w|=\max\{\|u\|,\|v\|\}.
\end{equation}
We observe that at $x\in U$
\begin{itemize}
\item if $w\in\FF_x^b$, then
  $|w|=\|v\|\le\|(u,v)\|=\|w\|\le\|u\|+\|v\|\le (1+b)\|v\|=(1+b)|w|$;
\item if $w\in\EE_x^a$, then
  $|w|=\|u\|\le\|(u,v)\|=\|w\|\le\|u\|+\|v\|\le (1+a)\|u\|=(1+a)|w|$.
\end{itemize}

Now for $x\in U$ such that $f_0(x)\in U$,
decomposing vectors in the $E$ and $F$ directions
$w=(u_0,v_0)\in\FF_x^b$ we have that $Df_0(x)\cdot w=(u_1,v_1)\in
\FF^{\lambda_0 b}_{f_0(x)}$ and also
\[
\|v_1\|\ge \|(Df_0\mid F(x))^{-1}\|^{-1}\cdot\|v_0\|\ge\|v_0\|
\qand
\|u_1\|\le \|Df_0\mid E(x)\|\cdot\|u_0\|\le \|u_0\|.
\]
Hence
\[
\frac{|Df_0(x)\cdot w|}{|w|} = \frac{\|v_1\|}{\|v_0\|} \ge
\|(Df_0\mid F(x))^{-1}\|^{-1}
\qand
\frac{\|Df_0(x)\cdot w\|}{\|w\|}
\ge \frac{\|(Df_0\mid F(x))^{-1}\|^{-1}}{1+b}.
\]
We observe that we can make the last expression as close to
$\|(Df_0\mid F(x))^{-1}\|^{-1}$ as we like by choosing $b$
very close to zero. Analogous calculations provide
\[
\frac{|Df_0^{-1}(x)\cdot w|}{|w|}   \ge
\|Df_0\mid E(f_0^{-1}(x))\|^{-1}
\qand
\frac{\|Df_0^{-1}(x)\cdot w\|}{\|w\|}
\ge
\frac{\|Df_0\mid E(f_0^{-1}(x))\|^{-1}}{1+a}
\]
for $x\in U$ such that $f_0^{-1}(x)\in U$ and $w\in\EE_x^a$.
Since the above calculations give approximately the same
bounds if we allow small perturbations in the factors
involved, then the same conclusion holds for other constants
$a^\prime,b^\prime$ perhaps closer to $0$ if we replace
$f_0$ by any sufficiently $C^1$-close map $g$.  We collect
this in the following lemma, which depends on the domination
assumption on the splitting, on the non-contractiveness
along $F$ and non-expansiveness along $E$, and also on the
isometric nature (specifically
property~\eqref{eq:invderivative})
of the perturbations we are considering.

\begin{lemma}
  \label{le:splitexpand}
  Let $f_0$ be a diffeomorphism admitting a dominated
  splitting $E\oplus F$ on a strictly forward invariant
  subset $U$ and $\Lambda=\cap_{n\ge0}\clos{f_0^n(U)}$. Let
  $\hat f:\cU\times M\to M$ be a family of isometric
  perturbations of $f_0$ as in Subsection
  \ref{sec:isom-rand-pert}. Then there exist
\begin{itemize}
\item angle bounds $a,b\in(0,1/2)$ defining stable
  $(\EE^a_x)_{x\in \clos(U)}$ and unstable $(\FF^b_x)_{x\in
    \clos(U)}$ cone fields in $T_{\clos(U)}M$;
\item a neighborhood $\cV$ of $0$ in $\cU$;
\item an open neighborhood $V$ of $\Lambda$ satisfying for
  every $v\in\cV$
\[
\clos(V)\subset U
,
\quad
\clos(f_v(V))\subset V
\qand
\clos(f_v^{-1}(V))\subset U;
\]
\end{itemize}
such that if $x\in V$ and
\begin{itemize}
\item $w\in\EE^a_x$, then
  $|Df_v^{-1}(x)\cdot w|\ge
  \|Df_0\mid E(f_v^{-1}(x))\|^{-1} \cdot |w|$;
\item $w\in\FF^b_x$, then
  $|Df_v(x)\cdot w|\ge
  \|(Df_0\mid F(x))^{-1}\|^{-1}\cdot |w|$.
\end{itemize}
\end{lemma}


\subsection{Uniform bound on the curvature of
  $E,F$-disks}
\label{sec:unif-bound-curv}

The ``curvature'' of the $E$- and $F$-disks defined at the
Introduction will be determined by the notion of H\"older
variation of the tangent bundle as follows.
Let us take $\de_0$ sufficiently small so that the
exponential map $\exp_x:B(x,\de_0)\to T_x M$ is a
diffeomorphism onto its image for all $x\in \clos(U_0)$, where
the distance in $M$ is induced by the Riemannian norm
$\|\cdot\|$. We write $V_x=B(x,\de_0)$ in what follows.  We
are going to identify $V_x$ through the local chart
$\exp_x^{-1}$ with the neighborhood $U_x=\exp_x (V_x)$ of
the origin in $T_x M$, and we also identify $x$ with the
origin in $T_x M$. In this way we get that $E(x)$ (resp.
$F(x)$) is contained in $\EE_y^a$ (resp. $\FF_y^b$) for all
$y\in U_x$, reducing $\de_0$ if needed, and the intersection
of $F(x)$ with $\EE_y^a$ (and the intersection of $E(x)$
with $\FF_y^b$) is the zero vector.

We write $\De$ also for the image of the respective
embedding for every $E$- or $F$-disk.  Hence if $\De$ is a
$E$-disk and $y=\De(w)$ for some $w\in\BB_E$, then the
tangent space of $\De$ at $y$ is the graph of a linear map
$A_x(y):T_x\De \to F(x)$ for $w\in\De^{-1}(V_x)$ (here
$T_x\De=D\De(x)(\RR^{d_E})$). The same happens locally for a
$F$-disk exchanging the roles of the bundles $E$ and $F$
above.

For $\ze\in(0,1)$ given by \eqref{eq:dominationclose} and
some $C>0$ we say that the \textit{tangent bundle of $\De$
  is $(C,\ze)$-H\"older} if
\begin{equation}
  \label{eq:holder-bundle}
  \| A_x(y) \| \le C \dist_\De(x,y)^\ze \quad\mbox{for all}\quad
y\in U_x\cap \De \qand x\in U,
\end{equation}
where $\dist_\De(x,y)$ is \textit{the distance along $\De$}
defined by the length of the shortest smooth curve from $x$
to $y$ inside $\De$ calculated with respect to the
Riemannian norm $\|\cdot\|$ induced on $TM$.

For a $E$- or $F$-disk $\De\subset U$ we define
\begin{equation}
  \label{eq:curvature-definition}
  \kappa(\De)=\inf\{C>0 : T\De \mbox{ is }(C,\ze)\mbox{-H\"older} \}.
\end{equation}
The proof of the following result can be easily adapted from
the arguments in~\cite[Subsection 2.1]{ABV00} with respect
to a single map $f_0$. The basic ingredients are the cone
invariance and dominated decomposition properties for $f_0$
that we have already extended for nearby diffeomorphisms
$g\in\cU$ with uniform bounds.

\begin{proposition}
\label{pr:bounded-curvature}
There is $C_1>0$ and a small neighborhood $X$ of $t_0$ such
that for every sequence $\omega\in\Omega=X^\NN$ with $X=\cU$

\begin{enumerate}
\item given a $F$-disk $\De\subset U$

  \begin{enumerate}
  \item there exists $n_1\in\NN$ such that
  $\kappa(f_\omega^n(\De))\le C_1$ for all $n\ge
  n_1$;
  \item if $\kappa(\De)\le C_1$ then
  $\kappa(f_\omega^n(\De))\le C_1$ for all
  $n\ge0$;
  \item in particular, if $\De$ is as in the previous item,
  then for every fixed $g=f_\omega$ with $\omega\in\Omega$
    $$
    J_n: f_\omega^n(\De) \ni x
    \mapsto \log | \det (Dg \mid T_x (f_\omega^n(\De)) |
    $$
    is $(L_1,\zeta)$-H\"older continuous with $L_1>0$
    depending only on $C_1$ and $f_0$, for every $n\ge1$.
  \end{enumerate}

\item for every $n\ge1$ and any given  $E$-disk $\De$ such that
  $(f_\omega^j)^{-1}(\De)\subset U$ for all $j=0,1,\dots,n$ and
  $\kappa(\De)\le C_1$, then
  \begin{enumerate}
  \item
  $\kappa\big( (f_\omega^n)^{-1}(\De) \big)\le
  C_1$ for all  $n\ge1$;
  \item for every $g=f_\omega$ with $\omega\in\Omega$ we have
    $$
    J_n: (f_\omega^n)^{-1}(\De) \ni x
    \mapsto \log | \det (Dg \mid
    T_x (f_\omega^n)^{-1}(\De)) |
    $$
    is $(L_1,\zeta)$-H\"older continuous with $L_1>0$
    depending only on $C_1$ and $f_0$.
  \end{enumerate}

\end{enumerate}
\end{proposition}

\begin{proof}
  See \cite[Proposition~2.2]{ABV00} and
  \cite[Corollary~2.4]{ABV00}.
\end{proof}

The bounds provided by
Proposition~\ref{pr:bounded-curvature} may be interpreted as
bounds on the curvature of either $E$-disks or $F$-disks,
since in the case $f_0\in\cU\subset\diff^2(M)$ we get $C_1$
as a bound on the curvature tensor of $\De$.


\subsection{Locally invariant sub-manifolds, expansion and contraction}
\label{sec:norm-invar-subm}

The domination assumption on $U_0$, the compactness and
$f_0$-invariance of $\Lambda$ together with properties
(1)-(2) from Theorem~\ref{th:limitSRB} ensure the existence
of families of $E$-disks ($C^{1+\zeta}$ center-stable manifolds)
$W^{cs}_{\de}(x)$ tangent to $E(x)$ at $x$ and $F$-disks
($C^{1+\zeta}$ center-unstable manifolds) $W^{cu}_\de (x)$ tangent
to $F(x)$ at $x$ which are locally invariant, for every
$x\in \Lambda$ and a small $\de>0$, as follows --- see
Hirsch-Pugh-Shub \cite{HPS77} for details.

There exist continuous families of embeddings
$\phi^{cs}:\Lambda\to\emb^{1+\zeta}(\BB_E,M)$ and
$\phi^{cu}:\Lambda\to\emb^{1+\zeta}(\BB_F,M)$, where
$\zeta\in(0,\alpha)$ is given by~\eqref{eq:dominationclose} and
$\emb^{1+\zeta}(\BB, M)$ is the space of $C^{1+\zeta}$ embeddings
from a ball $\BB$ in some Euclidean space to $M$, such that for all
$x\in\Lambda$
\begin{enumerate}
  \item $\phi^{cs}(x)$ is a $E$-disk and $T_x \phi^{cs}=E(x)$,
  $\phi^{cu}(x)$ is a $F$-disk and $T_x \phi^{cu}=F(x)$;
  \item writing $W^{cs}_{\de}(x)$ for
  $B(x,\de)\cap\phi^{cs}(x)(\BB_E)$ and $W^{cu}_\de (x)$ for
  $B(x,\de)\cap\phi^{cu}(x)(\BB_F)$ we have the local
  invariance properties: for every $\eta>0$ there exists $\de>0$
  such that for all $x\in\Lambda$
  \begin{enumerate}
  \item $f_0^{-1}(W^{cu}_\de(x))\subset
  W^{cu}_\eta(f_0^{-1}(x))$;
  \item $f_0(W^{cs}_\de(x))\subset
  W^{cs}_\eta(f_0(x))$.
  \end{enumerate}
\end{enumerate}


\subsubsection{Expansion/contraction of inner radius for $E/F$-disks}
\label{sec:InnerDiam}

Up to this point we have used some consequences of the
dominated decomposition assumption. Now we use assumptions
(1)-(2) of Theorem~\ref{th:limitSRB} to understand the
dynamical properties of the locally invariant sub-manifolds.

Given a smooth curve $\gamma:I\to M$ where $I=[0,1]$, we write
$\ell(\gamma)=\int_0^1\|\dot\gamma\|$ and
$L(\gamma)=\int_0^1|\dot\gamma|$ for the length of this
curve with respect to the norms $\|\cdot\|$ and $|\cdot|$.
Let
\[
\Gamma_E(\upsilon)=\{\gamma:I\to\BB_E: \gamma\mbox{ is
  smooth  and } 0<\|\dot\gamma\|\le\upsilon, \,
\gamma(0)=0, \, \gamma(1)\in\partial\BB_E\}
\]
and analogously for $\Gamma_F(\upsilon)$ with $\upsilon>0$.
We define the \emph{inner radius} of a $F$-disk $\De$ (with
respect to $|\cdot|$) to be
\[
R(\De)=\inf\{ L(\De\circ\gamma)\mid
\gamma\in\Gamma_F(\upsilon), \upsilon>0 \},
\]
and the \emph{inner diameter} of $\De$ to be
\[
\diam_\De(\De)=\sup\{ L(\De\circ\gamma)\mid
\gamma\in\Gamma_F(\upsilon), \upsilon>0 \},
\]
and likewise for $E$-disks. We note that $R=R(\De)\ge C
\dist(\De(0),\De(\partial \BB_F))>0$ where $C>0$ relates the
norms $\|\cdot\|$ and $|\cdot|$ and thus \emph{$R(\De)$ is a
minimum over $\Gamma_F(\upsilon)$ for some $\upsilon>0$}. For
fixing $\ep>0$ small we can find $\upsilon>0$ and
$\gamma\in\Gamma_F(\upsilon)$ such that $R\le
L(\De\circ\gamma)< R+\epsilon$, hence we may re-parametrize
$\gamma$ such that $\|(\De\circ\gamma)'\|$ is a constant in
$(C(R+\ep)^{-1},C R^{-1})$. Since $\Gamma_F(\upsilon)$ is a
compact family in the $C^1$ topology, we have that $R(\De)$
is assumed at some smooth curve.

Now we consider the family of $E$-disks having strictly
positive inner radius, bounded curvature and bounded inner diameter:
\begin{eqnarray*}
\cD_E(r,K,\de,k)
&=&
\{\De\in\emb^{1+\zeta}(\BB_E,M):
\quad
\De\mbox{  is a $E$-disk  },
\quad
\De(0)\in\clos{U_k},
\\
& &
\quad R(\De)\ge r,
\quad
\kappa(\De)\le K
\qand
\diam_\De(\De)\le\de\}
\end{eqnarray*}
for fixed $r,K,\de>0$ and $k\in\NN$, and analogously for $\cD_F(r,K,\de,k)$.

\begin{lemma}
  \label{le:compactfamilydisks}
  Given $r,K,\de>0$ and $k\in\NN$ the families
  $\cD_E(r,K,\de,k)$ and $\cD_F(r,K,\de,k)$ are compact in
  the $C^1$ topology of $\emb^{1+\zeta}(\BB_E, M)$ and
  $\emb^{1+\zeta}(\BB_F, M)$, respectively.
\end{lemma}

\begin{proof}
  We argue for $E$-disks only since the arguments for
  $F$-disks are the same.
  We note that $\cD_E(r,K,\de,k)$ defines a subset of bundle maps
  $D\De:T\BB_E\to TM, (x,v)\mapsto (\De(x),D\De(x)v)$. The
  bound on the ``curvature'' of the disks bounds the
  H\"older constant of $D\De(x)$ for $x\in\BB_E$. This H\"older
  control together with the bounded diameter condition
  $\int_0^1 |D\De(\gamma)\dot\gamma|\le\de$ ensures that
  $|D\De(x)|$ is equibounded on $\cD_E(r,K,\de,k)$. We also get
  that $\De(\BB_E)\subset B_1(\clos{U_0})$.

  Finally, the uniform bound on the H\"older constant of
  $D\De$ ensures that $D\De$ is a equicontinuous family for
  $\De\in\cD_E(r,K,\de,k)$. The proof finishes applying
  Ascoli-Arzela Theorem to $\{ D\De: \De\in\cD_E(r,K,\de,k)\}$ and
  noting that $\clos(U_0)$ is compact, any limit point must
  share the same inner radius and diameter bounds, and also
  that the cone families are continuous.
\end{proof}

From now on we fix $K=C_1$ from
Proposition~\ref{pr:bounded-curvature}, $k\in\NN$ big
enough, $\de>0$ small enough so that every $E$- and $F$-disk
in the above families be contained in $U_0$, and write
$\cD_E(r)$ and $\cD_F(r)$ for the families in
Lemma~\ref{le:compactfamilydisks}. Let $\lambda$ be
1-dimensional Lebesgue measure.

If we take $\De\in\cD_F(r)$ then $\cH(\De\cap F_1)<1$ by
assumption (2) from Theorem~\ref{th:limitSRB}.  Then
$\De\cap F_1$ is totally disconnected and \emph{curve free}
(see e.g.  \cite{falconer1990}), i.e. for any regular curve
$\gamma:I\to\De$ we have $\cH^1(\gamma(I)\cap F_1)=0$, where
$\cH^1$ is $1$-dimensional Hausdorff measure. Thus we must
have $\lambda(\gamma^{-1}(F_1))=0$. For otherwise
$\gamma(I)\cap F_1 = \gamma(\gamma^{-1}(F_1))$ and
$\cH^1(\gamma(I)\cap F_1)=\ell(\gamma\mid \gamma^{-1}(F_1))
= \int_{\gamma^{-1}(F_1)}\|\dot\gamma\|>0$, since $\gamma$
is a regular curve, a contradiction.

Lemma~\ref{le:splitexpand} and the fact that $\De$ is a
$F$-disk guarantee that
$L(g\circ\gamma) = \int_0^1
|Dg(\gamma(t)\cdot\dot\gamma(t)| >
\int_{\gamma^{-1}(\De\setminus F_1)}|\dot\gamma| +
\int_{\gamma^{-1}(F_1)} |\dot\gamma| = L(\gamma)$
for every $g\in\cU$ and smooth regular $\gamma:I\to\De$.
This is enough to show that $R(g(\De))>R(\De)$ for
$g\in\cU$, since $R(g(\De))$ is a minimum.
The compactness given by
Lemma~\ref{le:compactfamilydisks} assures that
there exists $\sigma_F=\sigma_F(r)>1$ such that
\begin{equation}
  \label{eq:rcu}
  R(g(\De))\ge\sigma_F \cdot R(\De)
\mbox{  for all  } \De\in\cD_F(r)\qand
g\in\clos{\cU},
\end{equation}
taking a smaller $\cU$ around $f_0$ if needed.
Clearly we can also get
$\sigma_E=\sigma_E(r)>1$ such that
\begin{equation}
  \label{eq:rcs}
  R(f_0^{-1}(\De))\ge\sigma_E \cdot R(\De)
\mbox{  for all  } \De\in\cD_E(r)\qand
g\in\clos{\cU},
\end{equation}
using the same arguments replacing $g$ by $g^{-1}$,
$F_1$ by $E_1$ and taking $\De\in\cD_E(r)$ above.

\begin{remark}
\label{rmk:localinv}
These estimates on the inner radius enable us to improve on
the local invariance properties from
Subsection~\ref{sec:norm-invar-subm} as follows: for every
$x\in\Lambda$ and $\de>0$ small enough there exists
$k=k(x,\de)\ge1$ satisfying
\begin{enumerate}
  \item $f_0^k(W^{cs}_\de(x))\subset
  W^{cs}_\de(f_0^k(x))$;
  \item $f_0^{-k}(W^{cu}_\de(x))\subset
  W^{cu}_\de(f_0^{-k}(x))$.
  \end{enumerate}
  Indeed for any given $\de>0$ we may find $\eta>0$ such
  that $R\big(f_0(W^{cs}_\de(x))\big)=\eta$ and so
\[
R(W^{cs}_\de(x))=R\big( f_0^{-1}( f_0(W^{cs}_\de(x)) ) \big)
\ge \sigma_E(\eta)\cdot R\big( f_0(W^{cs}_\de(x)) \big),
\]
thus $R\big( f_0(W^{cs}_\de(x)) \big)\le
\sigma_E(\eta)^{-1}\cdot R(W^{cs}_\de(x))<
R(W^{cs}_\de(x))$. This shows that for any given $\upsilon>0$
there must be an integer $k\ge1$ such that $R\big(
f_0^k(W^{cs}_\de(x)) \big)<\upsilon$, and analogously for the
center-unstable disks.
\end{remark}


\begin{remark}
\label{rmk:WinLambda}
In particular after Remark~\ref{rmk:localinv} we ensure that
$W^{cu}_\de(x)\subset \Lambda$ for every $x\in \Lambda$ and
$\de>0$ small enough.  For there is a constant $C>0$ (see
Subsection~\ref{sec:expand-contr-cone}) such that, if $y\in
W^{cu}_\de(x)$, then for every $\eta>0$ we have
$R(W^{cu}_\de(f_0^{-n}(x)))\le \eta$ and
$\dist(f_0^{-n}(x),f_0^{-n}(y))\le C\eta$ for big enough
$n\ge0$. Since $\Lambda$ is $f_0$-invariant we obtain
$\dist(f_0^{-n}(y),\Lambda)<C\eta$ or $y\in f_0^n(U_0)$ for
big $n\ge 0$ if $\eta$ is small enough. Hence $y\in\Lambda$.
Moreover this ensures that $W^{cu}_\de(x)$ is tangent to $F$
at every point.
\end{remark}


\subsection{A Local Product Structure for $\Lambda$}
\label{sec:Local-Prod-Struct}

The continuity in the $C^{1+\zeta}$ topology of $\phi^{cs}$ defined in
Subsection~\ref{sec:norm-invar-subm} and the inclusion
$W_\de^{cu}(x)\subset\Lambda$ obtained in
Remark~\ref{rmk:WinLambda} guarantee that for an open
neighborhood $V_0$ of $0$ in $\BB_F$ such that
$W_\de^{cu}(x)=\phi^{cu}(x)(V_0), \, x\in\Lambda$
\[
\psi_x:V_0\times \BB_E \to M, \quad
(y,z)\mapsto \phi^{cs}\big( \phi^{cu}(x)(y) \big) (z)
\]
is a $C^{1+\zeta}$ map for all $x\in\Lambda$. Moreover
\begin{itemize}
\item $D_1\psi_x(y,0)= D\big( \phi^{cu}(x) \big)
  (y):\RR^{d_F}\to F(y)$ is an isomorphism, since
  $\psi_x(y,0)=\phi^{cu}(x)(y)$ for all $y\in V_0$ and by
  definition of $\phi^{cu}$;
\item $D_2\psi_x(y,0):\RR^{d_E}\to E(y)$ is an isomorphism,
  by definition of $\phi^{cs}$.
\end{itemize}
Hence $|\det D\psi_x(y,0)|$ is bounded away from zero for
$y\in\Lambda$, because both the angle between $E(y)$ and
$F(y)$ (by domination), and $|\det D_2\psi_x(y,0)|$ are
bounded from below away from zero for $y\in\Lambda$ (by
compactness). Also we note that $\psi_x$ is just the
restriction to $W_\de^{cu}(x)\times \BB_E$ of a map
$\psi:\Lambda\times\BB_E\to M$. This shows that $\psi$ is a
local diffeomorphism from a neighborhood of $\Lambda\times
0$ in $\Lambda\times\BB_E$ to a neighborhood of $\Lambda$ in
$M$.  Since $\psi\mid(\Lambda\times 0)\equiv
\textrm{Id}\mid\Lambda$ we may choose a neighborhood $V_1$
of $0$ in $\BB_E$ so that $\psi_0=\psi\mid(\Lambda\times V_1)$ is a
diffeomorphism onto its image, which we write $W_0$.

\begin{remark}
\label{rmk:equalaverages}
In addition, following the arguments of Remark
\ref{rmk:WinLambda} we get
$\dist(f_0^{-n}(y),f_0^{-n}(x))\to0$ and
$\dist(f_0^{n}(z),f_0^{n}(x))\to0$ when $n\to+\infty$ for
all $x\in\Lambda, y\in W^{cu}_\de(x)$ and $z\in
W^{cs}_\de(x)$. In particular this shows that \emph{forward time
averages along center-stable disks and backward time averages along
center-unstable disks are constant}.
\end{remark}

The special neighborhood $W_0$ of $\Lambda$ together with
Remark~\ref{rmk:equalaverages} shows that the stable set
$W^s(\Lambda)=\{z\in M: \lim_{n\to+\infty}\dist(f_0^{n}(z),\Lambda)=0\}$ of
$\Lambda$ coincides with the union of the stable sets of
each point of $\Lambda$:  $W^s(\Lambda)=\cup_{x\in\Lambda} W^s(x)$.

\begin{lemma}
  \label{le:physicalcondition}
  There exist constants $h_0,h_1>0$ such that for any
  $f_0$-invariant ergodic probability measure $\mu$
  supported in $\Lambda$ and every $F$-disk
  $\De\subset\Lambda$, then $m(B(\mu))\ge h_0\cdot
  m_\De(B(\mu))$, where $m_\De$ is the Lebesgue measure
  induced by $m$ along the sub-manifold $\De$. In addition,
  if $\De=\De\cap B(\mu), m_\De-\bmod 0$, then there is a
  ball of radius $\ge h_1 \cdot R(\De)$ contained in
  $B(\mu)$ Lebesgue modulo zero and intersecting $\Lambda$.
\end{lemma}

\begin{proof}
  We have $B(\mu)\supset \cup_{y\in\De\cap B(\mu)}
  \phi^{cs}(y)(\BB_E)\supset \cup_{y\in\De\cap B(\mu)}
  W_\de^{cs}(y)$ by Remark~\ref{rmk:equalaverages} and
  definition of center-stable manifolds.  We note that both
  the angle between the tangent space to $x\in\De$ and
  $E(x)$, and the inner radius of the center-stable leaves
  are bounded from below away from zero over $\Lambda$ by
  $\beta_0$ and $r_0$ respectively. Hence the Lebesgue
  measure of $B_0=\cup_{y\in\De\cap B(\mu)} W_\de^{cs}(y)$
  is bounded from below by $h_0\cdot m_\De(B(\mu))$, where
  $h_0>0$ depends only on $\beta_0$ and $r_0$.  Thus if
  $\De\subset B(\mu), m_\De-\bmod0$, then $B_0$ contains a
  ball of radius bounded from below by $h_1\cdot R(\De)$
  dependent on $\beta_0$, on $r_0$ and on the curvature of
  $F$-disks, all uniform over $\Lambda$.  Clearly
  $B_0\cap\Lambda\neq\emptyset$.
\end{proof}


\subsubsection{Disks as graphs}
\label{sec:disks-as-graphs}

We can apply the results from
Subsection~\ref{sec:norm-invar-subm} to any sequence of maps
in $\cU$ using the invariance of $\hat\Lambda$ (see
Subsection~\ref{sec:RandomInvSet}) and the ``local product
structure'' from the previous discussion.

Let $K=C_1$ be as fixed in Subsection~\ref{sec:InnerDiam}.
Let $k\in\NN$ be big enough, $\de>0$ and $\ep_0$ small
enough be fixed so that every disk in $\cD_E(r), \cD_F(r)$
centered at $\hat\Lambda\subset \clos(U_k)\subset W_0$ be
contained in $W_0$, where $\hat\Lambda$ was defined in
Subsection~\ref{sec:RandomInvSet} for $\cU$ small enough
(corresponding to $\ep_0>0$ very small).
We consider the family of $E$- and $F$-disks which are
local graphs as follows
\begin{eqnarray*}
\cG_E(s)
&=&
\{\De\in\cD_E(r): r>0\mbox{ and for } x\in\Lambda\mbox{ such
that } \De(0)\in W^{cs}_\de (x),
\mbox{ there is } \phi:V\to V_0
\\
& & \qquad\qquad\quad
\mbox{ with }
B(0,s)\subset V\subset V_1\subset\BB_E\qand
\graph(\phi)\subset\psi_x^{-1}(\De)
\},
\end{eqnarray*}
and likewise for $\cG_F(s)$ with $s>0$, exchanging the roles
of $E, F, V_1$ and $V_0$.  We note that \emph{since cones are
complementary }(i.e. any $d_E$-subspace of $\EE_x$ together
with any $d_F$-subspace of $\FF_x$ span $T_x M$, $x\in U_0$)
then\emph{ every $E$- or $F$-disk is a local graph for some $s>0$}.
Let also $\de_0=\sup\{\diam\psi_x(V_0\times
V_1):x\in\Lambda\}>0$, which is finite by compactness.

\begin{lemma}
  \label{lem:bigleaves}
  Let $0<r\ll\min\{\de,\de_0\}$, $\De\in\cD_E(r)$ and
  $\hat\De\in\cD_F(r)$ be given. If
  $\De(0),\hat\De(0)\in\hat\Lambda$ and
  $\De\in\cG_E(s),\hat\De\in\cG_F(s)$ for some $s>0$, then
$(f_\omega)^{-1}(\De)\in\cG_E(\sigma_E(s)\cdot
  s), f_\omega(\hat\De)\in\cG_F(\sigma_F(s)\cdot s)$
for all $\omega\in\Omega$.

\end{lemma}

\begin{proof}
  It is obvious that $(f_\omega)^{-1}(\De)$ is a $E$-disk
  and that $f_\omega(\hat\De)$ is a $F$-disk after
  \eqref{eq:rcu} and \eqref{eq:rcs}. Moreover
  $(f_0)^{-1}(\De)\in\cG_E(\sigma_E(s)\cdot s),
  f_0(\hat\De)\in\cG_F(\sigma_F(s)\cdot s)$ by the local
  expression of $f_0$ on the ``local product coordinates''
  provided by $\psi$. The expansion on the inner radius of
  the domains of the graphs is a consequence of the fact
  that a ball in $W^{cs}_\de(x)$ is a $E$-disk and any ball
  in $W^{cu}_\de(x)$ is a $F$-disk, $x\in\Lambda$. The
  conclusion for $f_\omega$ and any $\omega\in\Omega$ holds
  since $f_\omega$ is taken $C^1$-close to $f_0$.
\end{proof}


\section{Equilibrium states and physical measures}
\label{sec:EquiStates}

Here we characterize the equilibrium states
$\mu$ for $f_0$ with respect to the potential $-\vfi(x)$
where $\vfi(x)=\log|\det Df_0\mid F(x)|$, as
in~\eqref{eq:3}.
We start by observing that, in the setting of
Theorem~\ref{th.limitequilibrium}, given any $f_0$-invariant
measure $\mu$ the sum $\chi^+(x)$ of the positive Lyapunov
exponents of $\mu$-a.e. point $x$ equals
\begin{equation}
    \label{eq:sumLyap}
\chi^+(x) = \lim_{n\to+\infty}\frac1n\log|\det Df^n_0\mid F(x)|
\end{equation}
by the Multiplicative Ergodic Theorem~\cite{Os68}. Indeed by
condition (1) of Theorem~\ref{th:limitSRB} every Lyapunov
exponent along the $E$ direction is non-positive and every
Lyapunov exponent along the $F$ direction is non-negative.

\begin{theorem}
\label{thm:equistates}
Let $f_0:M\to M$ be a $C^{1+\alpha}$ diffeomorphism
admitting a strictly forward invariant open set $U$ with a
dominated splitting satisfying items (1)-(2) of
Theorem~\ref{th:limitSRB}.  Then every equilibrium state
$\mu$ with respect to $\vfi$ as in~\eqref{eq:3}, supported
in $U$, is a convex linear combination of
\begin{enumerate}
\item finitely many ergodic equilibrium states with positive
  entropy and which are physical probability measures, with
    \item ergodic equilibrium states having zero entropy
      whose support has constant unstable Jacobian equal to
      one, i.e., measures whose Lyapunov exponents are
      non-positive.
\end{enumerate}
Moreover if $f_0\mid\Lambda$ is transitive, then there is at
most one ergodic equilibrium state with positive entropy
whose basin covers $U_0$ Lebesgue almost everywhere.
\end{theorem}

\begin{proof}
We first show that we may assume  $\mu$  ergodic.

\begin{lemma}\label{lem.ergdecompequil} Almost every ergodic
component of an equilibrium state for $\vfi$ is
itself an equilibrium state for the same function.  \end{lemma}

\begin{proof} Let $\mu$ be an $f$-invariant measure
  satisfying~\eqref{eq:3}.  On the one hand, the Ergodic
  Decomposition Theorem (see e.g Mañé~\cite{Man87}) ensures
  that
\begin{equation}
  \label{eq:9}
\int \vfi \,\, d\mu = \int\!\! \int \vfi(x) \, d\mu_z(x) \,
d\mu(z)
\quad\mbox{and}\quad h_\mu(f)=\int
h_{\mu_z}(f) \, d\mu(z).
\end{equation}
On the other hand, Ruelle's inequality guarantees for a
$\mu$-generic $z$ that (recall~\eqref{eq:sumLyap})
\begin{equation}
  \label{eq:10}
h_{\mu_z}(f)\le\int \vfi \, d\mu_z.
\end{equation}
By~\eqref{eq:9} and~\eqref{eq:10}, and because $\mu$ is an
equilibrium state~\eqref{eq:3}, we conclude that we have
equality in~\eqref{eq:10} for $\mu$-almost every $z$.
\end{proof}

Now let $\mu$ be an ergodic equilibrium state for $\vfi$
supported in $U$. Thus $\supp(\mu)\subset \Lambda$.  Now we
have two possibilities.
\begin{description}
\item[$h_\mu(f_0)>0$] According to the characterization of
  measures satisfying the Entropy Formula \cite{LY85}, $\mu$
  must be an $SRB$ \emph{measure}, i.e., \emph{$\mu$ admits
    a disintegration into conditional measures along
    unstable manifolds which are absolutely continuous with
    respect to the volume measure naturally induced on these
    sub-manifolds of $M$}.

\item[$h_\mu(f_0)=0$] Since $\vfi\ge0$ on $\Lambda$ by
  condition (1) of Theorem~\ref{th:limitSRB}, the
  equality~\eqref{eq:3} shows that $\vfi=0$ for $\mu$-a.e.
  $x$. Hence $\chi^+=0, \, \mu$-a.e. and $\mu$ has no
  expansion.
\end{description}

The Ergodic Decomposition then ensures that every
equilibrium state will be a convex linear combination of the
two types of measures described above. The latter
possibility corresponds to item (2) in the statement of
Theorem~\ref{thm:equistates}. The former case with positive
entropy needs more detail.

\begin{remark}
\label{rmk:CY-SRB}
Up until now we have shown that ergodic equilibrium states
for $-\vfi$ are either measures with no expansion or $SRB$
measures. This is exactly the same conclusion that
Cowieson-Young get \cite{CoYo2004} in a more general
setting.
\end{remark}

The Entropy Formula~\eqref{eq:3} and the assumption
$h_\mu(f_0)>0$ ensure that there are positive Lyapunov
exponents for $\mu$. Hence there exist \emph{Pesin's smooth
  ($C^{1+\alpha}$) unstable manifolds} $W^u(x)$ through
$\mu$-a.e. point $x$.  Moreover, as already mentioned, the
disintegration $\mu^u_x$ of $\mu$ along these unstable
manifolds $W^u(x)$ is absolutely continuous with respect to
the Lebesgue measure $m^u_x$ induced by the volume form of
$M$ restricted to $W^u(x)$, for $\mu$-a.e. $x$.

We claim that $\mu(F_1)=0$. For otherwise there would be
some component with $\mu_x^u(F_1)>0$ which implies
$m_x^u(F_1)>0$, and so
$\cH(F_1)\ge\cH(F_1\cap W^u(x))\ge\dim(W^u(x))\ge1$, a
contradiction, since $W^u(x)$ is a $F$-disk.

This means that $\int \log\|(Df_0\mid F)^{-1}\|\, d\mu <0$.
Hence \emph{the Lyapunov exponents of $\mu$ along every
  direction in $F$ are strictly positive}. Thus $\dim
W^u(x)=\dim F = d_F$ for  $\mu$-generic $x$, and $\mu$ is
a \emph{Gibbs state along the center-unstable direction}
$F$.  These manifolds are asymptotically backward
exponentially contracted by $Df_0$, see~\cite{Pe76}, hence
$W^u(x)\subset\Lambda$, since $\Lambda$ is a topological
attractor (see the arguments in Remark~\ref{rmk:WinLambda}).

Fixing a $\mu$-generic $x$, since $W^u(x)$ is an $F$-disk
Lemma~\ref{lem:bigleaves} ensures that we may assume
$R(W^u(x))\ge\rho$ for some $\rho>0$ dependent only on
$\dist(\Lambda,M\setminus U_0)$.  V\'asquez shows
\cite{vasquez2006} that the support of any Gibbs cu-state
such as $\mu$ contains \emph{entire unstable leaves
  $W^u(y)$} for $\mu$-a.e. $y$, so we may also assume that
$W^u(x)\cap B(\mu)$ has full Lebesgue measure in $W^u(x)$.
Using Lemma~\ref{le:physicalcondition} we get that $B(\mu)$
contains Lebesgue modulo zero a ball of radius uniformly
bounded from below by $h_1\cdot\rho>0$.

We have shown that each ergodic equilibrium state $\mu$
having positive entropy must be a physical measure. Since
the ergodic basins of distinct physical measures are
disjoint and have volume uniformly bounded from below away
from zero, there are at most finitely many such measures.
This concludes the proof of items (1) and (2) of
Theorem~\ref{thm:equistates}.

Let $f_0\mid \Lambda$ have a dense orbit and let us suppose
that there two distinct equilibrium states $\mu_1, \mu_2$
with positive entropy. Then by the previous discussion there
are two balls $B_1, B_2$ contained in the ergodic basins
$B(\mu_1)$ and $B(\mu_2)$ Lebesgue modulo zero,
respectively, and intersecting $\Lambda$. Since
$f_0\mid\Lambda$ is a transitive diffeomorphism and a
regular map, there exists $k\ge1$ such that $m(f^k(B_1)\cap
B_2)>0$. Thus $\mu_1=\mu_2$ and \emph{transitiveness of
  $\Lambda$ is enough to ensure there is only one
  equilibrium state with positive entropy}.

Let $\mu$ be the unique equilibrium state with
positive entropy and let us take an open set
$B=\psi_x(W^u(x)\times V_1)$ contained in $B(\mu)$ Lebesgue
modulo zero, where $x\in\Lambda\cap B(\mu)$ --- see
Subsection~\ref{sec:Local-Prod-Struct} for the definition of
$\psi_x$ on a $F$-disk such as $W^u(x)\subset\Lambda$. As
already explained, we may assume that $W^u(x)\cap B(\mu)$
has full Lebesgue measure along $W^u(x)$. We set
$\de=R(W^u(x))>0$.

\begin{figure}[htbp]
  \centering
  \includegraphics[width=12cm,height=3cm]{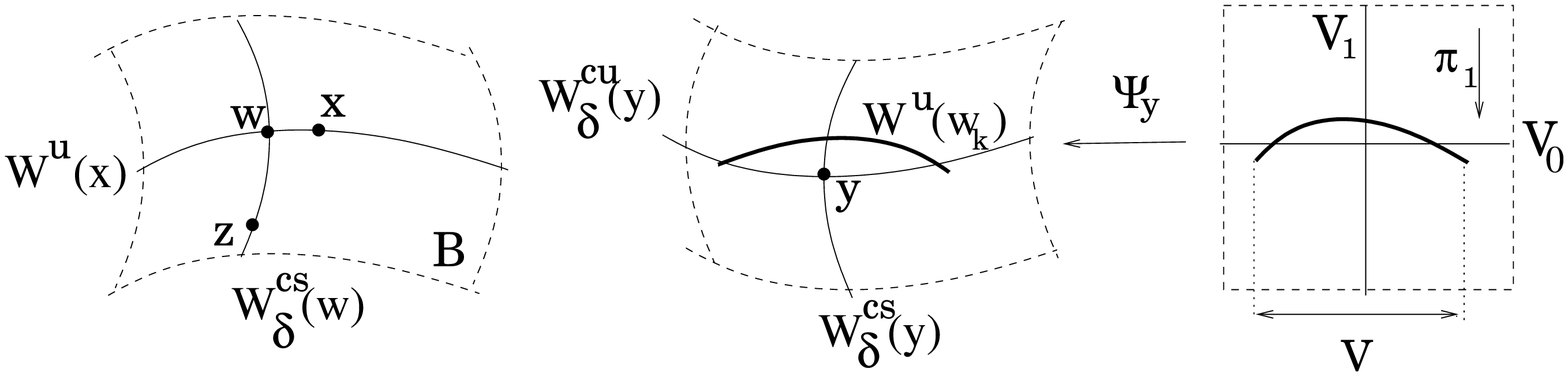}
  \caption{The ``local product structure'' neighborhoods $B$ and
$W^{cu}_\de(y)\times V_1$.}
  \label{figEPL}
\end{figure}

We can take $z\in B\cap\Lambda$ whose forward $f_0$-orbit is
dense in $\Lambda$.  We may take $z$ as close to $x$ as we
like and there are points $(w,u)\in W^u(x)\times V_1$ such
that $z=\psi_x(w,u)$. Then $z\in W^{cs}_\de(w)$ and
$W^{u}(w) = W^{u}(x)$, see Figure~\ref{figEPL}. Then by
Remark~\ref{rmk:equalaverages} the forward $f_0$-orbit of
$w$ is also dense in $\Lambda$. We note that
$f_0^n(W^{u}(w))\supset W^{u}(f_0^n(w))$ and also
$R(f_0^n(W^{u}(w))\ge\de$ for $n\ge1$ by
Lemma~\ref{lem:bigleaves}.

Arguing by contradiction, we suppose that $m(W_0\setminus
B(\mu))>0$, where $W_0$ was defined in
Subsection~\ref{sec:Local-Prod-Struct}.  Then there exists
$y\in\Lambda$ such that both $Z=(\psi_y\mid
(W^{cu}_\de(y)\times V_1))^{-1}(W_0\setminus B(\mu))$ and
$\pi_1(Z)$ have positive Lebesgue measure, where
$\pi_1:W^{cu}_\de(y)\times V_1\to W^{cu}_\de(y)$ is the
projection onto the first factor. Moreover we can choose $y$
so that it is a Lebesgue density point of $\pi_1(Z)$.

Let $n_k$ be a sequence such that $w_k=f^{n_k}_0(w)\to y$
when $k\to+\infty$. Since $W^u(w_k)$ is a $F$-disk,
$\psi_y^{-1}(W^u(w_k))\subset W^{cu}_\de(y)\times V_1$ is
the graph of a map from an open neighborhood $V$ of $y$ in
$W^{cu}_\de(y)$ to $V_1$, for big enough $k$. Then
$V\cap\pi_1(Z)$ has positive Lebesgue measure and so, after
Remark~\ref{rmk:equalaverages}, the Lebesgue measure of
$W^u(w_k)\cap(W_0\setminus B(\mu))$ is also positive. But
this implies that $W^u(x)\cap(W_0\setminus B(\mu))$ also has
positive Lebesgue measure, contradicting the choice of $x$.

This shows that $B(\mu)$ has full Lebesgue measure in $W_0$
and hence in $U_0$, as in the statement of
Theorem~\ref{thm:equistates}.
\end{proof}


\section{Zero-noise limits are equilibrium measures}
\label{sec:zero-noise-limits}

Here we prove Theorem~\ref{th.limitequilibrium}.  Let
$f_0:M\to M$, $\hat f: X\to C^{1+\alpha}(M,M), t\mapsto
f_t$, $f_{t_0}\equiv f$ for fixed $t_0\in X$, and
$(\th_\ep)_{\ep>0}$ be a family of probability measures on
$X$ such that $(\hat f,(\th_\ep)_{\ep>0})$ is a
non-degenerate isometric random perturbation of $f_0$, as in
Subsection~\ref{sec:isom-rand-pert}.

The main idea is to find a \emph{fixed random generating
  partition} for the system $(\hat f,\th_\ep)$ for
\emph{every small $\ep>0$} and use the absolute continuity
of the stationary measure $\mu^\ep$, together with the
conditions on the splitting to obtain a
semi-continuity property for entropy on zero-noise limits.

\begin{theorem}
  \label{thm.semicontentropy}
Let us assume that there exists a finite partition $\xi$ of $M$
(Lebesgue modulo zero) which is generating
for random orbits, for every small enough $\ep>0$.

Let $\mu^0$ be a weak$^*$ accumulation point of
$(\mu^\ep)_{\ep>0}$ when $\ep\to0$. If $\mu^{\ep_j}\to\mu^0$
for some $\ep_j\to 0$ when $j\to\infty$, then
$
\limsup_{j\to\infty} h_{\mu^{\ep_j}}(\hat f, \th_{\ep_j})
\le h_{\mu^0}(f_0,\xi).
$
\end{theorem}

\begin{remark}
\label{rmk:CYsimilar}
Recently Cowieson-Young obtained~\cite{CoYo2004} a
similar semi-continuity property without assuming the
existence of a uniform generating partition but using either
a local entropy condition or that the maps $\hat f$ involved
be of class $C^\infty$.
\end{remark}

The absolute continuity of $\mu^\ep$, the conditions on the
splitting for $f_0$ and the isometric
perturbations permit us to use a random version of the
Entropy Formula

\begin{theorem}
  \label{thm.randomPesin}
  If an ergodic stationary measure $\mu^\ep$ for a
  isometric random perturbation $(\hat f,\th_\ep)$ of
  $f_0$, in the setting of Theorem~\ref{th.limitequilibrium},
  is absolutely continuous for any given $\ep>0$, then
\[
h_{\mu^\ep}(\hat f, \th_\ep) =
\int \log|\det Df_0\mid F(x)| \, d\mu^\ep(x).
\]
Moreover if condition (3) of the statement of
Theorem~\ref{th:limitSRB} also holds, then in addition to
the above there exists $c>0$ such that
$h_{\mu^\ep}(\hat f, \th_\ep) \ge c$
for all $\ep>0$ small enough.
\end{theorem}

Putting Theorems~\ref{thm.semicontentropy}
and~\ref{thm.randomPesin} together shows that
$h_{\mu^0}(f_0)\ge \int \log|\det Df_0(x)| \, d\mu^0(x)$,
since $\th_\ep\to\delta_{t_0}$ in the weak$^*$ topology when
$\ep\to0$, by the assumptions on the support of $\th_\ep$ in
Subsection~\ref{sec:non-deg-pert}. Since the reverse
inequality holds in general (that is Ruelle's
inequality~\cite{Ru78}) we get the first statement of
Theorem~\ref{th.limitequilibrium}.  To conclude the proof we
just have to recall Theorem~\ref{thm:equistates} from
Section~\ref{sec:EquiStates}, which provides the second part
of the statement of Theorem~\ref{th.limitequilibrium}.


\subsection{Random Entropy Formula}
\label{sec:rand-pesins-form}

Now we explain how to obtain Theorem~\ref{thm.randomPesin}.

Let $\ep>0$ be fixed in what follows.  The Lyapunov
exponents $\lim_{n\to\infty} n^{-1}\log\|Df_\omega^n(x)\cdot
v\| $ exist for $\th^\ep\times\mu^\ep$-almost every
$(\omega,x)$ and every $v\in T_x M\setminus\{0\}$, by
Oseledets result \cite{Os68} adapted to this setting, see
e.g. \cite{arnold-l-1998}. At every given point $(\omega,x)$
there are at most $d=\dim(M)$ possible distinct values for
the above limit, the \emph{Lyapunov exponents} at
$(\omega,x)$. We write $\chi^+(\omega, x)$ for the \emph{sum
  of the positive Lyapunov exponents} at $x$. Lyapunov
exponents are $F$-invariant by definition, so
$\chi^+(\omega,x)=\chi^+(x)$ for
$\th^\ep\times\mu^\ep$-almost every $(\omega,x)$ (a
consequence of $\th^\ep$ being a product measure and
$\sigma$-ergodic, see~\cite[Corollary I.1.1]{LQ95}) and $\chi^+$ is
constant almost everywhere if $\mu^\ep$ is ergodic.

The Entropy Formula for random maps is the
content of the following result.
\begin{theorem}
  \label{thm:randomEntropyFormula}
Let a random perturbation $(\hat f,\th_\ep)$ of a
diffeomorphisms $f_0$ be given and assume that the
stationary measure $\mu^\ep$ is such that $\log|\det
Df_t(x)|\in L^1(\Omega\times M,\th_\ep\times\mu^\ep)$. If
$\mu^\ep$ is absolutely continuous with respect to Lebesgue
measure on $M$, then
\begin{equation}
  \label{eq:EF}
  h_{\mu^\ep}(\hat f,\th_\ep)=\int \chi^+\, d\mu^\ep.
\end{equation}
\end{theorem}

\begin{proof}
  See \cite{LeYo88} and \cite[Chpt. IV]{LQ95}.
\end{proof}

Now
since the random perturbations are isometric we have
\[
\frac1n\log\|(Df_\omega^n\mid F(x))^{-1}\|\le
\frac1n\sum_{j=0}^{n-1}
\log \|(Df_{\omega_{j+1}}\mid F(f^j_\omega(x)))^{-1}\|
\to \int \log\|(Df_0\mid F(x))^{-1}\| \,d\mu^\ep(x)
\]
when $n\to\infty$ for $\th^\ep\times\mu^\ep$-a.e.
$(\omega,x)$ by the Ergodic Theorem, if $\mu^\ep$ is ergodic.
By the assumptions on $f_0$ and $E\oplus F$ this
ensures that the Lyapunov exponents in the directions of $F$
are non-negative. In the same way we get for
$\th^\ep\times\mu^\ep$-a.e. $(\omega,x)$
\[
\limsup_{n\to\infty}
\frac1n\log\|Df_\omega^n\mid E(x)\|\le
 \int \log\|Df_0\mid E(x)\| \,d\mu^\ep(x)\le0,
\]
and so every Lyapunov exponent in the directions of $E$ is
non-positive. Since $E$ and $F$ together span $T_{U_0}M$,
according to the Multiplicative Ergodic Theorem
(Oseledets~\cite{Os68}) the sum $\chi^+$ of the positive
Lyapunov exponents (with multiplicities) equals the
following limit $\th^\ep\times\mu^\ep$-almost everywhere
\[
\chi^+(x)= \lim_{n\to\infty}\frac1n\log|\det Df_\omega^n\mid
F(x)| = \int\log |\det Df_0\mid F(x)|\,d\mu^\ep(x)\ge0.
\]
The identity above follows from the Ergodic Theorem, if
$\mu^\ep$ is ergodic, since the value of the limit is
$F$-invariant, thus constant.

Finally since $\mu^\ep$ is absolutely continuous for random
isometric perturbations, the formula in
Theorem~\ref{thm:randomEntropyFormula} gives the first part
of the statement of Theorem~\ref{thm.randomPesin}.

\begin{remark}
\label{rmk:uniformentropy}
The argument above together with conditions (1)-(3) from
Theorem~\ref{th:limitSRB} ensure that there exists $c_0>0$
satisfying $h_{\mu^\ep}(\hat f,\th_\ep)\ge c_0$ for every
small enough $\ep>0$. In fact, condition (3) ensures that
$|\det Df_0\mid F(x)|>1$ for all $x\in\Lambda$.  Hence there
is $c_0>0$ such that $\log|\det Df_0\mid F(x)|\ge c_0$ for
every $x$ in a neighborhood $U_k$ as in
Subsection~\ref{sec:RandomInvSet}, for some fixed big
$k\ge1$.
\end{remark}

Finally, as shown in Subsection~\ref{sec:RandomInvSet}, for
any given $k\ge1$ there is $\ep_0>0$ for which the random
invariant set $\hat\Lambda=\hat\Lambda_{\ep}$ is contained
in $U_k$ for all $\ep\in(0,\ep_0)$. Then
$\supp\mu^{\ep}\subset\hat\Lambda_\ep$ will be in the
setting of Remark~\ref{rmk:uniformentropy} above if
condition (3) of Theorem~\ref{th:limitSRB} holds in addition
to conditions (1) and (2). This completes the proof of
Theorem~\ref{thm.randomPesin}.


\subsection{Uniform random generating partition}
\label{sec:entr-with-gener}

Here we construct the uniform random generating partition
assumed in Theorem~\ref{thm.semicontentropy}.  In what
follows we fix a weak$^*$ accumulation point $\mu^0$ of
$\mu^\ep$ when $\ep\to 0$: there exist $\ep_j\to0$ when
$j\to\infty$ such that $\mu=\lim_{j\to+\infty} \mu^{\ep_j}$.


Let us take a finite cover $\{B(x_i,\rho_0/4),
i=1,\dots,\ell\}$ of $\clos(U_k)$ by $\rho_0/4$-balls, where
$\rho_0\in(0,\min\{ \de_0,\dist(M\setminus U_0, U_k)\} )$ for some $k\ge1$
such that $\supp(\mu^{\ep_j})\subset U_k\subset W_0$ for all
$j\ge1$. Recall from Subsection~\ref{sec:Local-Prod-Struct}
that $W_0$ is a ``local product structure'' neighborhood of
$\Lambda$ and note that we can choose $k$ as big as we like,
if we let $j$ be big enough.

Now since $\mu^0$ is a probability measure, we may assume
that $\mu^0(\partial \xi)=0$, for otherwise we can replace
each ball by $B(x_i,\gamma\rho_0/4)$, for some
$\gamma\in(1,3/2)$ and for all $i=1, \dots, k$.  We set
$\xi$ to be the finest partition of $M$ obtained through all
possible intersections of these balls:
$\xi=\{B(x_1,\gamma\rho_0/4), M\setminus
B(x_1,\gamma\rho_0/4)\}\vee\dots\vee
\{B(x_\ell,\gamma\rho_0/4), M\setminus
B(x_\ell,\gamma\rho_0/4)\}$. In what follows we let
$\rho=\gamma\rho_0/4\in(0,3\rho_0/4)$.

\begin{remark}
\label{rmk.0boundary}
The partition $\xi$ is such that all atoms of
$\vee_{i=-n}^{n} (f_\omega^i)^{-1} \xi$ have boundary
(which is a union of pieces of boundaries of open balls)
with zero Lebesgue measure, for all $n\ge1$ and every
$\omega\in\hat\Omega$. Moreover since $\mu^0$ is $f_0$-invariant
and $\mu^0(\partial \xi)=0$, then $\mu^0(\vee_{i=0}^{n-1}
f^{-j} \xi)=0$ for all $n\ge1$.
\end{remark}

\begin{lemma}\label{le.generatingpart}
For each $\omega\in\hat\Omega$ we have
$\diam\big(\vee_{i=-n}^{n} f_{\omega}^i (\xi)\big)\to0$
when $n\to+\infty$.
\end{lemma}

\begin{proof}
  Let $n\ge1$, $\omega\in\hat\Omega$, $x_0\in\clos(U_k)$ and
  $y_0\in (\vee_{i=-n}^{n} f_{\omega}^i (\xi))(x_0)$ with
  $y_0\neq x_0$. We write $x_k=f^k_\omega(x_0)$ and likewise
  for $y_k$, $|k|\le n$.

  Let us suppose that there exists a $E$-disk
  $\De\in\cG_E(s)$ centered at $x_0=\De(0)$ such that
  $y_0\in\De$ and $0<s_0<\dist(y_0,x_0)<\rho$.  Then by
  Lemma~\ref{lem:bigleaves} we see that since
  $\dist(y_i,x_i)\le\rho$ for $i=-1,\dots,-n$ we have
  $\dist(x_0,y_0)\le\sigma_E(s_0)^{-n}\cdot \rho$. If this
  holds for arbitrarily big values of $n\ge1$, then the
  statement of the lemma is proved.

  We now show that the assumption above is always true.
  Since $x_0,y_0\in\hat\Lambda$ we know that
  $x_n,y_n\in\hat\Lambda\subset W_0$. By definition we have
  $\dist(x_0,y_0), \dist(y_n,x_n)<\rho<\de_0$. Hence there
  exist $w_0,w_n\in\Lambda$ such that both
  $x_0,y_0\in\psi_{w_0}(V_0\times V_1)$ and $x_n\in
  W^{cs}_\de(w_n)$, and also $\De_n=W^{cs}_\de(w)$ is a
  $E$-disk and a graph through $x_n$, i.e.
  $\De_n\in\cG_E(\de)$.  Then applying
  Lemma~\ref{lem:bigleaves} several times we get
  $\De_0=(f^n_\omega)^{-1}(\De_n)\cap\psi_{w_0}(V_0\times
  V_1)\in\cG_E(\de)$.

\begin{figure}[htbp]
  \centering
  \includegraphics[width=9cm,height=3cm]{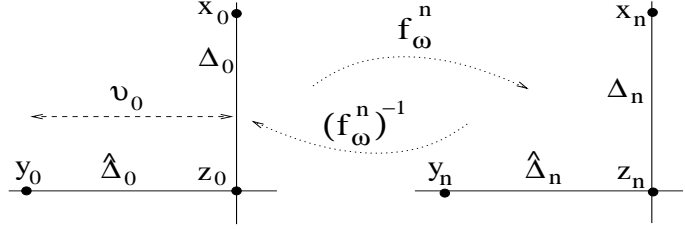}
  \caption{The construction of $\De_0, \hat\De_0, \De_n$ and
  $\hat\De_n$.}
  \label{fig:refinamento}
\end{figure}

  Now there exists $(y_0^u,y_0^s)\in V_0\times V_1$ such
  that $y_0=\psi_{w_0}(y_0^u,y_0^s)$. Let
  $\hat\De_0=\psi_{w_0}(V_0\times\{y_0^s\})\in\cG_F(\de)$ be
  a $F$-disk through $y_0$. Then we get a $F$-disk
  $\hat\De_n=f_\omega^n(\hat\De_0)\in\cG_F(\de)$ which is a
  graph through $y_n$. Since both
  $\De_n,\hat\De_n\subset\psi_{w_n}(V_0\times V_1)$ are
  graphs we know there exists a unique intersection $z_n$
  and thus there is
  $z_0=(f_\omega^n)^{-1}(z_n)=\De_0\cap\hat\De_0$, see
  Figure~\ref{fig:refinamento}.

  Let $\ups_0=\dist_{\hat\De_0}(y_0,z_0)$. We note that if
  $\ups_0=0$, then $y_0\in\De_0$ and we can proceed as in the
  beginning. Hence we assume $\ups_0>0$ and get $\rho> \dist(y_n,x_n) \ge
  \dist(y_n,z_n)-\dist(z_n,x_n) \ge \dist(y_n,z_n) - \rho$
  i.e. $\dist(y_n,z_n)<2\rho$. But by construction and
  applying Lemma~\ref{lem:bigleaves}
\[
\ups_0\le \dist_{\hat\De_0}(y_0,z_0) \le
\sigma_F(\ups_0)^{-n}\cdot \dist_{\hat\De_n}(y_n,z_n)
\le\sigma_F(\ups_0)^{-n}\cdot K_0\cdot \dist(y_n,z_n) \le
2\rho K_0 \cdot \sigma_F(\ups_0)^{-n},
\]
where $K_0$ is a constant relating distances in $M$ with
distances along $F$-disks and depending of the curvature
$\kappa(\hat\De_n)$, which is globally bounded, see
Subsection~\ref{sec:unif-bound-curv}.

This shows that $\ups_0$ can be made as small as we please.
Then for $\ups_0>0$ small enough there exists
$\De\in\cG_E(\de)$ with $x_0,y_0\in\De$, e.g. take the image
by $\psi_{w_0}$ of any  $d_E$-plane intersected with $V_0\times V_1\subset
\RR^d$ through  $(y_0^u,y_0^s), (0,y_0^s)$ and
$\psi_{w_0}^{-1}(x_0)$. Thus we can always reduce to the
first case above. The proof is complete.
 \end{proof}

 Lemma~\ref{le.generatingpart} implies that $\xi$ is a
 random generating partition Lebesgue modulo zero, hence
 $\mu^\ep$ modulo zero for all $\ep>0$, as in the statement
 of the Random Kolmogorov-Sinai Theorem~\ref{thm.KSrandom}.
 We conclude that $h_{\mu^{\ep_k}}((\hat
 f,\theta_{\ep_k}),\xi)=h_{\mu^{\ep_k}}(\hat f,
 \theta_{\ep_k})$ for all $k\ge1$.


\subsection{Semi-continuity of entropy on zero-noise}
\label{sec:semic-entr-zero-1}

Now we start the proof of Theorem~\ref{thm.semicontentropy}.
We need to construct a sequence of partitions of
$\hat\Omega\times M$ according to the following result ---
see Subsection~\ref{sec:metr-entr-rand} for the definitions
of $\hat\Omega$ and entropy. For a partition $\cP$ of a
given space $Y$ and $y\in Y$ we denote by $\cP(y)$ the
element (atom) of $\cP$ containing $y$.  We set
$\omega_0=(\dots,t_0,t_0,t_0,\dots)\in\hat\Omega$ in what
follows.

\begin{lemma}\label{lem.seqnpartitions} There exists a
sequence of measurable partitions $(\hat\cB_\ell)_{\ell\ge1}$ of
$\hat\Omega$ such that
 \begin{enumerate}
 \item $\omega_0\in\inter \hat\cB_\ell(\omega_0)$ for all
   $\ell\ge1$;
 \item $\hat\cB_\ell \nearrow \hat\cB$, $\hat\th^{\ep_j}
   \bmod 0$ for all $j\ge1$ when $n\to\infty$;
 \item $\lim_{n \to \infty} H_{\rho} (\xi \mid \hat\cB_n)
 = H_{\rho} (\xi \mid \hat\cB)$ for every measurable finite
 partition $\xi$ and any $G$-invariant probability measure $\rho$.
 \end{enumerate}
\end{lemma}

\begin{proof} For the first two items we let $\cC_n$ be a finite
  $\hat\th_{\ep_j}\bmod 0$ partition of $X$ such that
  $t_0\in\inter\cC_n(t_0)$ with $\diam \cC_n\to0$ when
  $n\to\infty$, for any fixed $j\ge1$. Example: take a cover
  $(B(t,1/n))_{t\in X}$ of $X$ by $1/n$-balls and take a
  sub-cover $U_1,\dots,U_l$ of $X\setminus B(t_0,2/n)$
  together with $U_0=B(t_0,3/n)$; then let
  $\cC_n=\{U_0, M\setminus U_0\}\vee\dots\vee
  \{U_l,M\setminus U_l\}$.

We observe that we may assume that the boundary of these
balls has null $\hat\th_{\ep_j}$-measure for all $j\ge1$, since
$(\hat\th_{\ep_j})_{j\ge1}$ is a denumerable family of
non-atomic probability measures on $X$ (see
Remark~\ref{rmk.noatoms}).  Now we set
\[
\hat\cB_n=
X^\NN\times\cC_n\times\stackrel{2n+1}{\dots}\times\cC_n\times
X^\NN\quad \mbox{for all }n\ge1,
\]
meaning that $\hat\cB_n$ is the family of all sets
containing points $\omega\in\hat\Omega$ such that
$\omega_i\in X$ for all $|i|>n$ and $\omega_i\in C_i$ for
some $C_i\in\cC_n, |i|\le n$.  Then since $\diam \cC_n\le
2/n$ for all $n\ge1$ we have $\diam \hat\cB_n\le 2/n$ and so
tends to zero when $n\to\infty$. Then $\hat\cB_n$ is an
increasing sequence of partitions and $\vee_{n\ge1}
\hat\cB_n$ generates the $\sigma$-algebra $\hat\cB,\,
\hat\th^{\ep_j} \bmod0$ (see e.g.  \cite[Lemma 3, Chpt.
2]{Bi65}) for all $j\ge1$. This proves items (1) and (2).
Item (3) is Theorem 12.1 of Billingsley~\cite{Bi65}.
\end{proof}

Now we deduce the right inequalities from known properties
of the conditional entropy. First we get from
Theorem~\ref{thm.randentropyFG} and~\cite[Thm. 0.5.3]{LQ95}
\begin{eqnarray*}
  h_{\mu^{\ep_j}}(\hat f,\theta_{\ep_j})
  &=&
  h_{\hat\mu^{\ep_j}}^{\hat\cB\times M}(G)
  = h_{ \hat\mu^{\ep_j}}^{\hat\cB\times M}(G,
  \hat\Omega\times\xi)
  \\
  &=&
  \inf \frac1n H_{\hat\mu^{\ep_j}}\left(
  \bigvee_{i=0}^{n-1} (G^i)^{-1}(\hat\Omega\times\xi) \mid
  \hat\cB\times M \right),
\end{eqnarray*}
where $\hat\Omega\times\xi=\{ \hat\Omega\times A: A\in\xi\}$.
Then for any given fixed $N\ge1$ and for every $\ell\ge1$
\begin{eqnarray*}
 h_{\mu^{\ep_j}}(\hat f,\theta_{\ep_j})
 &\le&
 \frac1N H_{\hat\mu^{\ep_j}}\left(
  \bigvee_{i=0}^{N-1} (G^i)^{-1}(\hat\Omega\times\xi) \mid
  \hat\cB\times M \right)
 \\
 &\le&
 \frac1N H_{\hat\mu^{\ep_j}}\left(
  \bigvee_{i=0}^{N-1} (G^i)^{-1}(\hat\Omega\times\xi) \mid
  \hat\cB_\ell\times M \right)
\end{eqnarray*}
because $\hat\cB_\ell\times M\subset\hat\cB\times M$.  Now
we fix $N$ and $\ell$, let $j\to \infty$ and note that since
\[
\mu^0(\partial\xi)=0=\de_{\omega_0}(\partial \hat\cB_m)
\quad\mbox{then}\quad
(\de_{\omega_0}\times\mu^0)(\partial(B_i\times\xi_l))=0
\]
for all $B_i\in\hat\cB_m$ and $\xi_l\in\xi$, where
$\de_{\omega_0}$ is the point mass concentrated at
$\omega_0$.  By weak$^*$ convergence
$\th^{\ep_l}\to\de_{\omega_0}$ and $\mu^{\ep_l}\to\mu^0$ we
get
$\hat\mu^{\ep_l}\to\hat\mu^0=\delta_{\omega_0}\times\mu^0$
when $l\to\infty$, see Lemma~\ref{le:weakcontinuous}.
Hence
\begin{equation}
  \label{eq:5}
  \limsup_{j\to\infty} h_{\mu^{\ep_j}}(\hat
  f,\theta_{\ep_k})
  \le
  \frac1N H_{\de_{\omega_0}\times\mu^0}\left(
  \bigvee_{i=0}^{N-1} (G^i)^{-1}(\hat\Omega\times\xi) \mid
  \hat\cB_\ell\times M \right)
  =\frac1N H_{\mu^0} \big(\bigvee_{i=0}^{N-1} f_0^{-i}\xi \big).
\end{equation}
Here it is easy to see that the middle conditional entropy
of~\eqref{eq:5} (involving only finite partitions) equals
$N^{-1}\sum_i \mu^0(P_i)\log\mu^0(P_i)$, where
$P_i=\xi_{0}\cap f^{-1}\xi_{1}\cap\dots\cap f^{-(N-1)}
\xi_{N-1}$ ranges over every sequence of possible atoms
$\xi_{0},\dots,\xi_{N-1}\in \xi$.

Finally, since $N$ was an arbitrary integer,
Theorem~\ref{thm.semicontentropy} follows from the
inequality in~\eqref{eq:5}. As already explained, this
completes the proof of Theorem~\ref{th.limitequilibrium}.


\section{Stochastic stability}
\label{sec:stochstability}

Here we prove Theorem~\ref{th:limitSRB}. Let $f_0:M\to M$ be
as in the statement of Theorem~\ref{th.limitequilibrium} and
let $\mu$ be an equilibrium state for $-\vfi$, as
in~\eqref{eq:3} (recall the definition of $\vfi$ in
Section~\ref{sec:EquiStates}) obtained using the
construction described in Section~\ref{sec:Isom-pert-maps}
through non-degenerate random isometric perturbations.

Condition (3) in the statement of Theorem~\ref{th:limitSRB}
ensures that the only possibility for the ergodic
decomposition of $\mu$ is the one given by item (2a) in
statement of Theorem~\ref{th.limitequilibrium}. In fact,
after Remark~\ref{rmk:uniformentropy}, every weak$^*$
accumulation point $\mu$ of $\mu^\ep$ when $\ep\to0$ will be
not only an equilibrium state for $-\vfi$, as shown in
Section~\ref{sec:zero-noise-limits}, but will also have
strictly positive entropy $h_\mu(f_0)\ge c>0$, after the
statement of Theorem~\ref{thm.randomPesin}.
Hence combining the statements in
Section~\ref{sec:zero-noise-limits} with
Theorem~\ref{th.limitequilibrium} we see that \emph{every
  weak$^*$ accumulation point $\mu$ of $\mu^\ep$ when
  $\ep\to0$ is a finite convex linear combination of the
  ergodic equilibrium states for $-\vfi$, which are physical
  measures}.

This shows that the family of equilibrium states for $-\vfi$
in the setting of Theorem~\ref{th:limitSRB} is
stochastically stable.

In addition, if $f_0$ is transitive, then there is only one
equilibrium state $\mu$ for $-\vfi$ which is ergodic and
whose basin covers $U_0$ Lebesgue almost everywhere, by the
last part of the statement of Theorem~\ref{thm:equistates}.
Then \emph{every weak$^*$ accumulation point of $\mu^\ep$
  when $\ep\to0$ necessarily equals $\mu$}. This finishes
the proof of Theorem~\ref{th:limitSRB}.



\bibliographystyle{plain}


\begin{thebibliography}{10}

\bibitem{ABV00}
J.~F. Alves, C.~Bonatti, and M.~Viana.
\newblock {SRB} measures for partially hyperbolic systems whose central
  direction is mostly expanding.
\newblock {\em Invent. Math.}, 140(2):351--398, 2000.

\bibitem{Ze03}
Jos{\'e} Alves.
\newblock {\em Statistical analysis of non-uniformly expanding dynamical
  systems}.
\newblock Publica\c c\~oes Matem\'aticas do IMPA. [IMPA Mathematical
  Publications]. Instituto de Matem\'atica Pura e Aplicada (IMPA), Rio de
  Janeiro, 2003.
\newblock 24$\sp {\rm o}$ Col\'oquio Brasileiro de Matem\'atica. [24th
  Brazilian Mathematics Colloquium].

\bibitem{AA03}
Jose~F. Alves and Vitor Araujo.
\newblock Random perturbations of nonuniformly expanding maps.
\newblock {\em Ast\'erisque}, 286:25--62, 2003.

\bibitem{Ar00}
V.~Ara\'ujo.
\newblock Attractors and time averages for random maps.
\newblock {\em Annales de l'Inst. Henri Poincar\'e - Analyse Non-lin\'eaire},
  17:307--369, 2000.

\bibitem{ArTah}
V.~Ara\'ujo and A.~Tahzibi.
\newblock Stochastic stability at the boundary of expanding maps.
\newblock {\em Nonlinearity}, 18:939--959, 2005.

\bibitem{arnold-l-1998}
Ludwig Arnold.
\newblock {\em Random dynamical systems}.
\newblock Springer-Verlag, Berlin, 1998.

\bibitem{BaV96}
V.~Baladi and M.~Viana.
\newblock Strong stochastic stability and rate of mixing for unimodal maps.
\newblock {\em Ann. Sci. \'Ecole Norm. Sup. (4)}, 29(4):483--517, 1996.

\bibitem{BeV2}
Michael Benedicks and Marcelo Viana.
\newblock Random perturbations and statistical properties of {H}\'enon-like
  maps.
\newblock {\em Ann. Inst. H. Poincar\'e Anal. Non Lin\'eaire}, 23(5):713--752,
  2006.

\bibitem{Bi65}
P.~Billingsley.
\newblock {\em Ergodic theory and information}.
\newblock J. Wiley \& Sons, New York, 1965.

\bibitem{BoV00}
C.~Bonatti and M.~Viana.
\newblock S{R}{B} measures for partially hyperbolic systems whose central
  direction is mostly contracting.
\newblock {\em Israel J. Math.}, 115:157--193, 2000.

\bibitem{BR75}
R.~Bowen and D.~Ruelle.
\newblock The ergodic theory of {A}xiom {A} flows.
\newblock {\em Invent. Math.}, 29:181--202, 1975.

\bibitem{brin-kifer1987}
M.~Brin and Yu. Kifer.
\newblock Dynamics of markov chains and stable manifolds for random
  diffeomorphisms.
\newblock {\em Ergodic Theory and Dynamical Systems}, 7:351--374, 1987.

\bibitem{Ca93}
M.~Carvalho.
\newblock {S}inai-{R}uelle-{B}owen measures for $n$-dimensional derived from
  {A}nosov diffeomorphisms.
\newblock {\em Ergod. Th. \& Dynam. Sys.}, 13:21--44, 1993.

\bibitem{CatEnr2001}
E~Catsigeras and E.~Enrich.
\newblock {S}{R}{B} measures of certain almost hyperbolic diffeomorphisms with
  a tangency.
\newblock {\em Discrete and Continuous Dynamical Systems}, 7(1):177--202, 2001.

\bibitem{CoYo2004}
W.~Cowieson and L.-S. Young.
\newblock {SRB} measures as zero-noise limits.
\newblock {\em Ergodic Theory and Dynamical Systems}, 25(4):1115--1138, 2005.

\bibitem{falconer1990}
Kenneth Falconer.
\newblock {\em Fractal Geometry: mathematical foundations and applications}.
\newblock John Wiley \& Sons, New York, USA, 1990.

\bibitem{FHY83}
A.~Fathi, M.-R. Herman, and J.-C. Yoccoz.
\newblock A proof of {P}esin's stable manifold theorem.
\newblock In {\em Geometric dynamics (Rio de Janeiro, 1981)}, volume 1007 of
  {\em Lecture Notes in Math.}, pages 177--215. Springer, Berlin, 1983.

\bibitem{guillemin-pollack1974}
V.~Guillemin and A.~Pollack.
\newblock {\em Differential Topology}.
\newblock Prentice Hall, New Jersey, 1974.

\bibitem{HPS77}
M.~Hirsch, C.~Pugh, and M.~Shub.
\newblock {\em Invariant manifolds}, volume 583 of {\em Lect. Notes in Math.}
\newblock Springer Verlag, New York, 1977.

\bibitem{hirsch1976}
Morris Hirsch.
\newblock {\em Differential Topology}.
\newblock Springer-Verlag, New-York, 1976.

\bibitem{hu2000}
H.~Hu.
\newblock Conditions for the existence of {SBR} measures for ``almost
  {A}nosov'' diffeomorphisms.
\newblock {\em Trans. Amer. Math. Soc.}, 352(5):2331--2367, 2000.

\bibitem{hu2001}
H.~Hu.
\newblock Statistical properties of some almost hyperbolic systems.
\newblock In {\em Smooth ergodic theory and its applications (Seattle, WA,
  1999)}, volume~69 of {\em Proc. Sympos. Pure Math.}, pages 367--384. Amer.
  Math. Soc., Providence, RI, 2001.

\bibitem{hu-young1995}
H.~Hu and L.-S. Young.
\newblock Nonexistence of {SBR} measures for some diffeomorphisms that are
  ``almost {A}nosov''.
\newblock {\em Ergodic Theory Dynam. Systems}, 15(1):67--76, 1995.

\bibitem{Ki86a}
Yu. Kifer.
\newblock General random perturbations of hyperbolic and expanding
  transformations.
\newblock {\em J. Analyse Math.}, 47:11--150, 1986.

\bibitem{Ki86}
Yuri Kifer.
\newblock {\em Ergodic theory of random transformations}, volume~10 of {\em
  Progress in Probability and Statistics}.
\newblock Birkh\"auser Boston Inc., Boston, MA, 1986.

\bibitem{Ki88}
Yuri Kifer.
\newblock {\em Random perturbations of dynamical systems}, volume~16 of {\em
  Progress in Probability and Statistics}.
\newblock Birkh\"auser Boston Inc., Boston, MA, 1988.

\bibitem{LY85}
F.~Ledrappier and L.-S. Young.
\newblock The metric entropy of diffeomorphisms {I}. characterization of
  measures satisfying {P}esin's entropy formula.
\newblock {\em Ann. of Math}, 122:509--539, 1985.

\bibitem{LeYo88}
F~Ledrappier and L.-S. Young.
\newblock Entropy formula for random transformations.
\newblock {\em Probab. Theory and Related Fields}, 80(2):217--240, 1988.

\bibitem{LQ95}
P.-D. Liu and M.~Qian.
\newblock {\em Smooth ergodic theory of random dynamical systems}, volume 1606
  of {\em Lect. Notes in Math.}
\newblock Springer Verlag, 1995.

\bibitem{Man87}
R.~Ma{{\~n}}\'e.
\newblock {\em Ergodic theory and differentiable dynamics}.
\newblock Springer Verlag, New York, 1987.

\bibitem{Nash1954}
John Nash.
\newblock {$C\sp 1$} isometric imbeddings.
\newblock {\em Ann. of Math. (2)}, 60:383--396, 1954.

\bibitem{Nash1956}
John Nash.
\newblock The imbedding problem for {R}iemannian manifolds.
\newblock {\em Ann. of Math. (2)}, 63:20--63, 1956.

\bibitem{Os68}
V.~I. Oseledets.
\newblock A multiplicative ergodic theorem: {L}yapunov characteristic numbers
  for dynamical systems.
\newblock {\em Trans. Moscow Math. Soc.}, 19:197--231, 1968.

\bibitem{Pe76}
Ya. Pesin.
\newblock Families of invariant manifolds corresponding to non-zero
  characteristic exponents.
\newblock {\em Math. USSR. Izv.}, 10:1261--1302, 1976.

\bibitem{Ru76}
D.~Ruelle.
\newblock A measure associated with {A}xiom {A} attractors.
\newblock {\em Amer. J. Math.}, 98:619--654, 1976.

\bibitem{Ru78}
David Ruelle.
\newblock {An inequality for the entropy of differentiable maps.}
\newblock {\em Bol. Soc. Bras. Mat.}, 9:83--87, 1978.

\bibitem{Si72}
Ya. Sinai.
\newblock Gibbs measures in ergodic theory.
\newblock {\em Russian Math. Surveys}, 27:21--69, 1972.

\bibitem{thaler1980}
M.~Thaler.
\newblock Estimates of the invariant densities of endomorphisms with
  indifferent fixed points.
\newblock {\em Israel Journal of Mathematics}, 37(4):303--314, 1980.

\bibitem{thaler1983}
M.~Thaler.
\newblock Transformations on $[0,1]$ with infinite invariant measures.
\newblock {\em Israel Journal of Mathematics}, 46(1-2):67--96, 1983.

\bibitem{Vi97b}
Marcelo Viana.
\newblock {\em Stochastic dynamics of deterministic systems}.
\newblock Publica\c c\~oes Matem\'aticas do IMPA. [IMPA Mathematical
  Publications]. Instituto de Matem\'atica Pura e Aplicada (IMPA), Rio de
  Janeiro, 1997.
\newblock 21$\sp {\rm o}$ Col\'oquio Brasileiro de Matem\'atica. [21th
  Brazilian Mathematics Colloquium].

\bibitem{vasquez2006}
Carlos Vásquez.
\newblock Statistical stability for diffeomorphisms with dominated splitting.
\newblock {\em Ergodic Theory and Dynamical Systems}, 27(1):253--283, 2007.

\bibitem{williams1970}
R.~F. Williams.
\newblock The "{DA}" maps of smale and structural stability.
\newblock {\em Proc. Symp. Pure Math., Amer. Math. Soc.}, 14:329--334, 1970.

\bibitem{Yo85}
L.-S. Young.
\newblock Stochastic stability of hyperbolic attractors.
\newblock {\em Ergod. Th. \& Dynam. Sys.}, 6:311--319, 1986.

\end{thebibliography}

\def\cprime{$'$}

\end{document}